\documentclass[11pt,a4paper,reqno]{article}
\usepackage[margin=1.8cm]{geometry}
\usepackage{amsfonts,amsthm,amsmath,amscd,amssymb}
\usepackage{indentfirst,enumitem,graphicx,epstopdf,subcaption}
\usepackage{pifont, dsfont,bm,mathtools,xcolor,empheq,float,hyperref,tikz-cd}
\usepackage{microtype}
\hypersetup{hidelinks,pdftitle={Recovering fractal boundary dimensions from Robin eigenvalues}}
\numberwithin{equation}{section}
\theoremstyle{plain}
\usepackage{booktabs}
\graphicspath{{plot/}}
\newtheorem{thm}{Theorem}[section]
\newtheorem{lemma}[thm]{Lemma}
\newtheorem{prop}[thm]{Proposition}
\newtheorem{cor}[thm]{Corollary}
\theoremstyle{remark}
\newtheorem{remark}[thm]{Remark}
\theoremstyle{definition}
\newtheorem{definition}[thm]{Definition}

\newcommand{\nm}{\noalign{\smallskip}}
\newcommand{\ds}{\displaystyle}
\begin{document}
\title{Recovering fractal boundary dimensions from Robin eigenvalues}
\date{}
\author{
    Habib Ammari\thanks{Department of Mathematics, ETH Zürich, Rämistrasse 101, 8092 Zürich, Switzerland and Hong Kong Institute for Advanced Study, City University of Hong Kong, Kowloon Tong, Hong Kong Special Administrative Region. \texttt{habib.ammari@math.ethz.ch}.} 
    \and
    Yat Tin Chow \thanks{Department of Mathematics, University of California, Riverside, Riverside, CA, USA. \texttt{yattinc@ucr.edu}.}
    \and 
    Fuqun Han \thanks{Department of Mathematics, City University of Hong Kong, Kowloon Tong, Hong Kong Special Administrative Region. \texttt{fuqun.han@cityu.edu.hk}  (corresponding author).} 
}
\maketitle

\begin{abstract}
We study Robin problems in which a small local boundary coefficient competes with the increasing surface measure of a finite-depth fractal approximation. For a fixed smooth domain, layer potentials give Green function and eigenpair expansions; contour residues and a finite-rank reduction yield a quantitative spectral-gap-dependent eigenvalue remainder. For varying Lipschitz domains, a uniform Poincar\'e inequality and a trace estimate with respect to normalized surface measure give a uniform ground-state expansion and a rank-one resolvent singularity as the total Robin boundary mass tends to zero. Under the boundary-resolution hypothesis and the condition on the approximation scales, ground-state measurements with known Robin coefficients recover the lower and upper Minkowski dimensions for sequences with bounded total Robin boundary mass. We describe Neumann, generalized Robin, and Dirichlet-type operator limits under the corresponding compact Mosco convergence hypotheses. Koch snowflake computations illustrate the reconstruction formulas, and an exact finite-stage area correction improves the adjacent-stage dimension estimator.
\end{abstract}


\bigskip


\section{Motivation and main results}\label{sec:introduction}

The detection of inaccessible corrosion is an important problem in non-destructive evaluation.  
Boundary measurements and modal data offer
complementary ways of probing changes that cannot be observed directly.
In particular, asymptotic methods have been used to identify internal corrosion from electrostatic measurements, vibration data, and ultrasound
measurements; see \cite{corrosion-num,corrosion-vibration,corrosion-wave}.
The present paper studies a different, multiscale aspect of this problem:
the interaction between an increasingly rough boundary and a small local Robin coefficient.

Corrosion morphology can exhibit irregular pit distributions and roughness over a range of scales.  Fractal and multifractal descriptors
have been used to quantify such features in corrosion images and in comparisons with electrochemical measurements
\cite{corrosion-ref1,corrosion-ref2,corrosion-ref3}.  Garc\'ia-Ochoa and
Corvo \cite{GarciaOchoaCorvo2010} compare fractal descriptors of copper
patina images and electrochemical-noise signals.  In the related setting of rough and porous surfaces, Risovi\'c and collaborators
\cite{RisovicEtAl2008} compare dimensions obtained from scanning electron
microscopy and electrochemical impedance spectroscopy.  Fractal porous
models also arise in simulations of corrosion-driven damage
\cite{corrosion-eth}.  We idealize this multiscale geometry by finite-depth
prefractal boundaries.  This does not assert that every corroded surface has an exact self-similar law, or that its electrochemical response is determined by a fractal dimension alone.

We consider the scalar Robin model with boundary condition
$\partial_\nu u+\varepsilon u=0$, where $\varepsilon>0$ is a real boundary exchange coefficient after nondimensionalization.  We refer to $\varepsilon$ as the Robin coefficient. Here, $\varepsilon$ represents local boundary exchange rather than electrochemical impedance.  Both the geometry and this effective coefficient may depend on corrosion time.  Our analysis is formulated for coupled sequences
$(\Omega_m,\varepsilon_m)$, so it accommodates either
$\delta=\delta(\varepsilon)$ or $\varepsilon=\varepsilon(\delta)$ when such a parametrization is appropriate.  The spectral interpretation is that of a scalar model; converting measured physical frequencies into its eigenvalues requires the corresponding wave speed and constitutive parameters to be calibrated.  A related physical motivation comes from Laplacian transport across irregular interfaces: Filoche and Sapoval \cite{FilocheSapoval2000} study transfer across random and deterministic prefractal geometries, while Grebenkov, Filoche, and Sapoval \cite{GrebenkovFilocheSapoval2006} formulate the interface response through the Dirichlet-to-Neumann operator.  
These transport models motivate the competition between surface geometry and local exchange.
 
For a fixed bounded smooth connected domain, the limit $\varepsilon\to0$ is a perturbation of the Neumann problem.  The zero Neumann eigenvalue is simple and its eigenfunction is constant. 
Thus, as $\varepsilon\to0$, the ground-state Robin eigenvalue tends to zero, while the Green operator evaluated at zero spectral parameter has a leading singular term of order $\varepsilon^{-1}$ associated with the constant Neumann eigenfunction. With the sign convention $\Delta_xG_\varepsilon(x,y)=\delta_y$, the leading terms are
\[
 \lambda_\varepsilon^0 =\varepsilon\frac{|\partial\Omega|}{|\Omega|} +\mathcal O(\varepsilon^2),
 \qquad
 G_\varepsilon(x,y) =-\frac1{\varepsilon|\partial\Omega|}+\mathcal O(1),
\]
where the second statement is understood away from the diagonal, or through the regular-part estimate \eqref{green-leading-expansion}. For a fixed simple positive Neumann eigenvalue, the first two Robin coefficients and the corresponding eigenfunction correction are expressed by boundary products of Neumann modes.  Appendix~\ref{sec:proofs} proves these formulas directly from layer-potential characteristic values, contour residues, and finite-rank factorizations in the framework of Ammari, Kang, and Lee \cite{surv153,surv235}.  
It also treats multiple eigenvalues and the zero mode, proves convergence of the coefficient series, and provides a spectral-gap-dependent remainder estimate.
For a broader account of the Robin spectral problem, we refer to Bucur, Freitas, and Kennedy \cite{bucur}.  These fixed-mode results are distinct from high-energy Robin--Neumann comparisons and spectral averages, for which we refer to \cite{rudnick2021,rupert2025}.
The associated boundary restriction questions are treated by Tataru \cite{tataru} and, with spectral weights adapted to Neumann boundary values, by Barnett, Hassell, and Tacy \cite{BarnettHassellTacy2018}.

The main difficulty for a varying boundary is that a fixed-domain remainder need not be uniform in the prefractal depth. 
If $\Gamma_m=\partial\Omega_m$, with surface measure $\sigma_m$ and trace operator $T_m$, the boundary contribution is
\begin{equation}\label{effective-mass-measure}
 \varepsilon_m\int_{\Gamma_m}|T_m u|^2\,d\sigma_m
       =\tau_m\int_{\Gamma_m}|T_m u|^2\,d\widehat\sigma_m,
 \qquad
 \tau_m:=\varepsilon_m|\Gamma_m|,
 \quad \widehat\sigma_m:=|\Gamma_m|^{-1}\sigma_m.
\end{equation}
The normalized measure $\widehat\sigma_m$ is a probability measure, and $\tau_m$ is the total Robin boundary mass.  Consequently, $\varepsilon_m\to0$ does not imply $\tau_m\to0$, because the growth of the boundary measure $|\Gamma_m|$ can offset the decrease in the local Robin coefficient. 
We seek estimates whose constants and range of validity are independent of the approximation level $m$. The normalized trace estimate provides this control when the expansion is expressed in terms of $\tau_m$. Simply substituting a varying boundary measure into a fixed-domain Taylor formula does not ensure that its remainder is controlled independently of $m$.
Extension and trace results on irregular domains provide a relevant background; see Jones \cite{Jones1981SobolevExtension} and Biegert \cite{Biegert2009Traces}.
The uniform weighted estimates of Capitanelli and Vivaldi \cite{CapitanelliVivaldi2014} concern Dirichlet problems on snowflake prefractal domains.  They illustrate the need for depth-independent control, but are distinct from the normalized Robin trace assumption used here.

Our first main result is a uniform estimate that requires only
volumes bounded above and away from zero and a common Poincar\'e-trace inequality with
respect to $\widehat\sigma_m$.  In the notation
$\Lambda_m(\varepsilon)=\lambda^0_{\varepsilon,\delta_m}$,
Lemma~\ref{lem:uniform-prefractal-ground-state} establishes
\[
 c_0\min\{\varepsilon|\Gamma_m|,1\}
       \leq\Lambda_m(\varepsilon)
       \leq\frac{\varepsilon|\Gamma_m|}{|\Omega_m|},
\]
and
\[
 \Lambda_m(\varepsilon)
       =\frac{\varepsilon|\Gamma_m|}{|\Omega_m|}
          \bigl(1+\mathcal O(\varepsilon|\Gamma_m|)\bigr)
       \quad\hbox{as }\varepsilon|\Gamma_m|\to0,
\]
with constants independent of $m$.  The corresponding resolvent estimate in Corollary~\ref{cor:subcritical-green-singularity} identifies the same constant-mode singularity uniformly in the prefractal depth.

Our second main result concerns the spectral reconstruction of boundary growth.  Assume that the approximating surface measure resolves the covering numbers of the limiting boundary $\Gamma$ at scale $\delta_m$, in the sense that $|\Gamma_m|\asymp\delta_m^{d-1}\mathcal N(\Gamma,\delta_m)$. If the approximation scales satisfy the logarithmic spacing condition \eqref{scale-sampling-condition} and $\sup_m\tau_m<\infty$, Theorem~\ref{thm:minkowski-recovery} recovers the lower and upper Minkowski dimensions from the logarithmic growth of $\Lambda_m(\varepsilon_m)/\varepsilon_m$.  Standard self-similar constructions provide concrete examples of the boundary-resolution hypothesis \cite{Hutchinson1981SelfSimilarity,Falconer2014FractalGeometry}.

The dimension reconstruction requires known Robin coefficients $\varepsilon_m$. When the domain volumes $|\Omega_m|$ are also known, subcritical eigenvalue measurements yield $|\Gamma_m|\sim|\Omega_m|\Lambda_m(\varepsilon_m)/\varepsilon_m$. If the volumes are known but the Robin coefficients are unknown, the same leading-order data determine only the total Robin boundary mass $\tau_m=\varepsilon_m|\Gamma_m|\sim|\Omega_m|\Lambda_m(\varepsilon_m)$, without separating the boundary measure from the Robin coefficient. The reconstruction uses measurements across approximation levels; it does not establish uniqueness for the joint recovery of an arbitrary shape and an unknown Robin coefficient from a single eigenvalue.

This distinction is also reflected in related inverse results.
Sincich \cite{Sincich2010} proves logarithmic stability for simultaneous recovery of an unknown corroded boundary and its impedance from two measurements on an accessible boundary.  
Santacesaria and Yachimura \cite{SantacesariaYachimura2020} recover a Robin coefficient in a mixed Dirichlet--Robin problem from the principal eigenvalue together with normal-derivative data of the principal eigenfunction on the Dirichlet part of the boundary.
Choulli, Metidji, and Soccorsi \cite{ChoulliMetidjiSoccorsi2024} study the recovery of an interior potential for a Robin Laplacian from boundary spectral data, consisting of eigenvalues and boundary traces of eigenfunctions.  
Joint shape and generalized-impedance reconstruction from far-field scattering data is studied by Bourgeois, Chaulet, and Haddar \cite{BourgeoisChauletHaddar2012}.  These works concern different unknowns or richer measurement sets; they do not replace the assumptions of known Robin coefficients and boundary resolution required by our multiscale eigenvalue result.

We also distinguish total-mass estimates from operator convergence. Generalized Robin problems on arbitrary domains are treated by Daners \cite{Daners2000} and, using boundary measures and capacity, by Arendt and Warma \cite{ArendtWarma2003}.  
Bass, Burdzy, and Chen \cite{BassBurdzyChen2008} investigate how fractal geometry affects mixed Robin problems through a probabilistic approach.  Neumann conditions on fractal boundaries are studied by Achdou and Tchou \cite{AchdouTchou2007}. 
For prefractal approximations, Capitanelli \cite{Capitanelli2010,Capitanelli2010Mixed} establishes Robin convergence results, including mixed Dirichlet--Robin problems, and Capitanelli and Vivaldi \cite{CapitanelliVivaldi2011} study insulating layers on Koch mixtures. 
The conductive transmission-layer homogenization of Lancia, Mosco, and Vivaldi \cite{LanciaMoscoVivaldi2008} is a related, but distinct, example of a limiting energy with a singular interface contribution.

We place the convergence of energy forms in the framework developed by Mosco \cite{Mosco1994}.  The uniform-domain results of Hinz, Rozanova-Pierrat, and Teplyaev \cite{HinzRozanovaTeplyaev2023} provide Poincar\'e inequalities and variational convergence for suitably converging domains and boundary measures. Under the measure-convergence and compact Mosco hypotheses specified in Section~\ref{sec:operator-regimes}, $\tau_m\to0$ gives a Neumann limit, whereas $\tau_m\to\tau_*\in(0,\infty)$ gives a generalized Robin limit. If $\tau_m\to\infty$, the limit is Dirichlet-type provided the required zero-trace and recovery conditions also hold.

The paper is organized as follows.  Section~\ref{sec:fixed-domain} collects the fixed-domain Green function and eigenpair formulas and their quantitative frequency limitations.  
Section~\ref{sec:prefractal} proves the uniform prefractal estimates and Minkowski-dimension recovery, and states the additional assumptions needed for the operator limits.
Section~\ref{sec:applications} discusses self-similar scaling and the Koch computations.  Its numerical analysis separates an exact finite-stage area correction from discretization and perturbation errors.  
Appendix~\ref{sec:proofs} provides the complete layer-potential proofs of the eigenpair expansions and the quantitative remainder.

\section{Fixed smooth domain}\label{sec:fixed-domain}

\subsection{Green function}

Let $\Omega\Subset\mathbb R^d$, $d=2,3$, be a bounded connected domain with boundary of class $C^{1,\alpha}$, $0<\alpha<1$.  Consider the Green function $G_{\varepsilon }$ defined by
\begin{equation} \label{def:green}
\begin{cases}
\Delta_x G_{\varepsilon }(x,y) = \delta_y \text{ in } \Omega, \\
\nm\ds
\frac{\partial G_{\varepsilon }(x,y)}{\partial \nu(x)}+ \varepsilon  G_{\varepsilon }(x,y) = 0 \quad \text{on } \partial \Omega,
\end{cases}
\end{equation}
for \(\varepsilon  > 0\). Integrating $\Delta G_\varepsilon =\delta_y$ and using the Robin boundary condition gives 
\begin{equation} \label{green:normalization}
   \int_{\partial \Omega} G_\varepsilon(x,y)\, d\sigma(x) = - \frac{1}{\varepsilon}.
\end{equation}

Next, we introduce the Neumann function $N$
\begin{equation}
\label{def:Neumann}
\begin{cases}
\Delta_x N(x,y) = \delta_y \quad \text{ in } \Omega, \\
\nm\ds
\frac{\partial N(x,y)}{\partial \nu(x)} = \frac{1}{|\partial \Omega|} \quad\text{ on } \partial \Omega, \\
\nm\ds
\int_{\partial \Omega} N(x,y) \, d\sigma(x) = 1.
\end{cases}
\end{equation}

The normalization \eqref{green:normalization} suggests subtracting the
constant singular term before solving for the harmonic correction.
We first introduce the layer-potential spaces used in this construction.

Let $S_\Omega$ denote the Laplace single-layer potential and its
boundary trace, and let $K_\Omega^*$ be the adjoint
Neumann--Poincar\'e operator. Its interior normal-derivative jump relation is
\[
\partial_\nu S_\Omega[\phi]\big|_-
=\left(-\frac12I+K_\Omega^*\right)[\phi].
\]

Let $H^*$ be the real space $H^{-1/2}(\partial\Omega)$ equipped with the inner product
\[
(\phi,\psi)_{H^*}
=
\begin{cases}
-(\phi,S_\Omega[\psi])_{H^{-1/2},H^{1/2}},&d=3,\\
-(\phi,\widetilde S_\Omega[\psi])_{H^{-1/2},H^{1/2}},&d=2.
\end{cases}
\]
In dimension two, the unmodified logarithmic single layer may have
a one-dimensional kernel. We therefore use the modified single layer
\[
\widetilde S_\Omega[\phi]
=S_\Omega\!\left[
\phi-\left(\int_{\partial\Omega}\phi\,d\sigma\right)\phi_e
\right]
-\left(\int_{\partial\Omega}\phi\,d\sigma\right)\chi_{\partial\Omega},
\]
where $\phi_e$ is the equilibrium density with total mass one.
This operator agrees with $S_\Omega$ on zero-total-mass densities;
see \cite[Section~2.4]{surv235}.

Recall that $K_\Omega^*:H^*\to H^*$ is compact and self-adjoint and let
\[
\{\phi_j,\kappa_j\}
\]
 denote eigenpairs where $\{\phi_j\}_{j\geq 0}$ is an orthonormal basis for $H^*$. We take $\kappa_0=1/2$; then $\kappa_j\to0$ as $j\to\infty$. Define
\[
H_0^*
:=
\left\{\phi\in H^*:
(\phi,\chi_{\partial\Omega})_{H^{-1/2},H^{1/2}}=0\right\}
=
\left\{\phi\in H^*:(\phi,\phi_0)_{H^*}=0\right\}.
\]
In both dimensions, set $\widetilde G_\varepsilon:=G_\varepsilon-N$ and $F_\varepsilon:=-|\partial\Omega|^{-1}\chi_{\partial\Omega}-\varepsilon N(\cdot,y)$ on $\partial\Omega$. To illustrate the role of the equilibrium mode, we first assume $d=3$; the calculation through \eqref{j>0} is restricted to this case.
Since $\widetilde G_\varepsilon$ is harmonic, we write
\begin{equation}\label{G:rep} \widetilde G_\varepsilon(x,y)=S_\Omega[\phi_\varepsilon](x). \end{equation}
The Robin boundary condition gives
\begin{equation}
      \label{eqn_phi} \left(-\frac12I+K_\Omega^*\right)[\phi_\varepsilon]+\varepsilon S_\Omega[\phi_\varepsilon]=F_\varepsilon \quad\text{on }\partial\Omega. 
\end{equation}
For the $H^*$-normalized eigenfunction $\phi_0$ associated with $\kappa_0=1/2$, we have
\begin{equation}\label{phi0}
       S_\Omega[\phi_0]=-\frac{1}{\int_{\partial\Omega}\phi_0\,d\sigma}\quad\text{on }\partial\Omega.
 \end{equation}

To solve \eqref{eqn_phi}, we first compute
\begin{align*}
(F_{\varepsilon}, \phi_{0})_{H^{*}} = &\int_{\partial \Omega} \phi_{0}(x) S_{\Omega} [\frac{1}{|\partial \Omega|} \mathcal{\chi}_{\partial \Omega} + \varepsilon N(\cdot,y)](x) \, d\sigma(x) \\
= & \int_{\partial \Omega} S_{\Omega} [\phi_0](x) (\frac{1}{|\partial\Omega|} \mathcal{\chi}_{\partial \Omega} + \varepsilon N(x,y)) \, d\sigma(x)\\
= & - \frac{(1+\varepsilon)}{\int_{\partial \Omega} \phi_0\, d\sigma}\,,
\end{align*}
where we have used \eqref{phi0}.

The spectral decomposition of $K_\Omega^*$ gives 
\[
K_{\Omega}^* [\phi] = \sum_{j=0}^{+\infty} \kappa_j (\phi, \phi_j)_{H^*} \phi_j.
\]

To solve \eqref{eqn_phi}, we expand the density in the eigenbasis:
\[
\phi_{\varepsilon } = \sum_{j=0}^{+\infty} (\phi_{\varepsilon}, \phi_{j})_{H^{*}} \phi_{j}
\]
and obtain
\[
\sum_{j=0}^{+ \infty} (-\frac{1}{2} + \kappa_j)(\phi_{\varepsilon }, \phi_{j})_{H^*} \phi_{j} + \varepsilon  S_{\Omega} [\phi_{\varepsilon }] = F_{\varepsilon }.
\]
Taking the $H^*$ inner product with $\phi_\ell$ gives
\[
(-\frac{1}{2} + \kappa_{\ell}) (\phi_{\varepsilon }, \phi_{\ell })_{H^{*}}  + \varepsilon  (S_{\Omega} [\phi_{\varepsilon }], \phi_{\ell })_{H^{*}}  = (F_{\varepsilon }, \phi_{\ell })_{H^{*}}. 
\]

Then,  for \( \phi_{0} \), we have
\begin{align} \label{j=0}
\varepsilon \bigg[- \frac{|\partial \Omega|}{(\int_{\partial \Omega} \phi_0\, d\sigma)^2} (\phi_{\varepsilon}, \phi_0)_{H^*} + \frac{1}{\int_{\partial \Omega} \phi_0\,d\sigma}  \sum_{j>0} (\phi_{\varepsilon}, \phi_j)_{H^*} \int_{\partial \Omega} \phi_j S_{\Omega}[\chi_{\partial \Omega}]\,d\sigma     \bigg] = - \frac{(1+\varepsilon)}{\int_{\partial \Omega} \phi_0\,d\sigma},
\end{align}
and, for  \( \phi_{j}, j>0\), 
\begin{align}\label{j>0}
    (- \frac{1}{2} + \kappa_j) (\phi_{\varepsilon}, \phi_j)_{H^*} + \varepsilon \sum_{\ell=0}^{+\infty} (\phi_{\varepsilon}, \phi_{\ell})_{H^*} (S_{\Omega}[\phi_{\ell}], \phi_j)_{H^*} = (F_\varepsilon, \phi_j)_{H^*}.
\end{align}
Equation~\eqref{j=0} isolates the equilibrium-mode coupling. By \eqref{green:normalization} and \eqref{def:Neumann}, $\int_{\partial\Omega}(G_\varepsilon-N)\,d\sigma=-(1+\varepsilon)/\varepsilon$. We therefore separate the constant $-(1+\varepsilon)/(\varepsilon|\partial\Omega|)$, leaving a harmonic correction with zero boundary mean. The following centered construction applies to both $d=2$ and $d=3$.

\begin{prop}[Green function for a small Robin coefficient]
\label{prop:green-fixed}
For every fixed $y\in\Omega$, define 
\begin{align} \label{def:f0}
F_\varepsilon^0(x,y) := \varepsilon \Big( \frac{1}{|\partial \Omega|} - N(x,y) \Big) \in H^*_0,
\end{align}
and 
\begin{align}
\label{G:approx-projected}
G_\varepsilon^{\mathrm{app}}(x,y) :={}&N(x,y)-\frac{1+\varepsilon} {\varepsilon|\partial\Omega|} +\sum_{j=1}^{+\infty}\frac{(F_\varepsilon^0,\phi_j)_{H^*}}{\kappa_j-\frac12}\left[S_\Omega[\phi_j](x)-\frac1{|\partial\Omega|}\int_{\partial\Omega}\phi_jS_\Omega[\mathcal\chi_{\partial\Omega}]\,d\sigma\right].
\end{align}
Then the (correction) series converges in $H^1(\Omega)$ and
\begin{equation}\label{green-remainder}
\|G_\varepsilon(\cdot,y)
-G_\varepsilon^{\mathrm{app}}(\cdot,y)\|_{H^1(\Omega)}
\leq C_y\varepsilon^2,
\qquad 0<\varepsilon\leq\varepsilon_0.
\end{equation}
The constant can be chosen uniformly when $y$ ranges over a compact subset
of $\Omega$.  Consequently,
\begin{equation}\label{green-leading-expansion}
\left\|G_\varepsilon(\cdot,y)-N(\cdot,y)
+\frac{1+\varepsilon}{\varepsilon|\partial\Omega|}\right\|_{H^1(\Omega)}
\leq C_y\varepsilon.
\end{equation} 
Boundary integrals containing an $H^{-1/2}$ density in the displays above are understood as $H^{-1/2}$--$H^{1/2}$ duality pairings.
\end{prop}

\begin{proof}
For $y \in \Omega$, we write
\[
G_\varepsilon(x,y)
=N(x,y)-\frac{1+\varepsilon}
{\varepsilon|\partial\Omega|}+R_\varepsilon(x,y).
\]
It follows from \eqref{green:normalization} that
$\int_{\partial\Omega}R_\varepsilon\,d\sigma=0$ and hence $R_{\varepsilon} \big|_{\partial \Omega} \in H^*_0$, while the boundary
conditions for $G_\varepsilon$ and $N$ give
\begin{equation}\label{eq:R-epsilon}
\begin{cases}
\Delta_x R_\varepsilon (x,y) =0&\text{in }\Omega,\\
\nm\ds
\frac{\partial R_\varepsilon (x,y)}{\partial \nu(x)} +\varepsilon R_\varepsilon(x,y)=F_\varepsilon^0(x,y) &\text{on }\partial\Omega,
\end{cases}
\end{equation}
where $F_\varepsilon^0$ is defined by \eqref{def:f0}.

Define $S_{\Omega,0}: H^* \rightarrow H^{1/2}(\partial \Omega)$ by
\begin{equation} \label{def:S0}
S_{\Omega,0}[\phi]:= S_\Omega[\phi]-\frac1{|\partial\Omega|}\int_{\partial\Omega}\phi S_\Omega[\mathcal\chi_{\partial\Omega}]\,d\sigma,
\end{equation}
and use the same notation for its extension in $\Omega$.  
Then there is a unique $\psi_\varepsilon\in H_0^*$ such that $R_\varepsilon= S_{\Omega,0}[\psi_\varepsilon]$.  The jump relation turns
\eqref{eq:R-epsilon} into
\begin{equation}\label{eq:psi-epsilon}
\left(-\frac12I+K_\Omega^*+\varepsilon  S_{\Omega,0}\right)[\psi_\varepsilon]=F_\varepsilon^0\quad\text{in }H_0^*.
\end{equation}
The operator $-\frac12I+K_\Omega^*$ is boundedly invertible on $H_0^*$. Since $F_\varepsilon^0=\mathcal O_{H^*}(\varepsilon)$, a Neumann-series argument in \eqref{eq:psi-epsilon} gives
\begin{equation}\label{psi-expansion}
\psi_\varepsilon =\sum_{j=1}^{+\infty} \frac{(F_\varepsilon^0,\phi_j)_{H^*}} {\kappa_j-\frac12}\phi_j +\mathcal O_{H^*}(\varepsilon^2).
\end{equation}
Applying $S_{\Omega,0}:H^*\to H^1(\Omega)$ and using \eqref{def:S0} proves \eqref{G:approx-projected} and \eqref{green-remainder}. 
\end{proof}

\begin{remark} We first note that 
$$
\int_{\partial \Omega} G_\varepsilon^{\mathrm{app}}(x,y)\, d\sigma(x) = - \frac{1}{\varepsilon}, 
$$
which is consistent with \eqref{green:normalization}.
Next, note also that if $\Omega$ is a disk or a ball, then $\mathcal\chi_{\partial\Omega}\in\operatorname{span}\{\phi_0\}$, or equivalently, $S_\Omega[\mathcal\chi_{\partial\Omega}]$ is constant on $\partial\Omega$. Hence, 
\begin{equation}\label{constant-density-condition}
(\mathcal\chi_{\partial\Omega},\phi_j)_{H^*}=0, \qquad j\geq1,
\end{equation}
which shows that $(F_\varepsilon^0,\phi_j)_{H^*}=(F_\varepsilon,\phi_j)_{H^*}$, and hence \eqref{G:approx-projected} becomes the unprojected formula
\begin{align}\label{G:approx}
G_\varepsilon^{\mathrm{app}}(x,y) ={}&N(x,y)-\frac{1+\varepsilon}{\varepsilon|\partial\Omega|}+\sum_{j=1}^{+\infty}\frac{(F_\varepsilon,\phi_j)_{H^*}}{\kappa_j-\frac12}\left[S_\Omega[\phi_j](x)-\frac1{|\partial\Omega|}\int_{\partial\Omega}\phi_jS_\Omega[\mathcal\chi_{\partial\Omega}]\,d\sigma\right],
\end{align}
with the error estimate \eqref{green-remainder}.  
\end{remark}

\subsection{Eigenvalue problem}

We let $0=\lambda_0^0<\lambda_0^1 \leq \lambda_0^2 \leq \dots$ be the eigenvalues of $-\Delta$ on $\Omega$ with Neumann boundary conditions on $\partial \Omega$. Let $u^j_{0}$ denote the normalized (real) eigenfunction associated with $\lambda_0^j$; that is, it satisfies $\| u^{j}_0\|_{L^2(\Omega)}=1$.

We consider the eigenvalue problem:
\begin{equation} \label{pb:eigenvalue}
\begin{cases}
- \Delta u^j_{\varepsilon} = \lambda_{\varepsilon}^{j} u^j_{\varepsilon} & \text{in } \Omega, \\
\nm\ds
\frac{\partial u^j_{\varepsilon}}{\partial \nu} + \varepsilon u^j_{\varepsilon} = 0 & \text{on } \partial \Omega.
\end{cases}
\end{equation}

For each $\varepsilon\geq0$, we choose the eigenfunctions $u_\varepsilon^j$ to form a real orthonormal basis of $L^2(\Omega)$.

Our aim is to provide asymptotic expansions of \(\lambda_{\varepsilon}^{j}\) as \(\varepsilon \to 0\) and  $j$ fixed. We recall first that the Robin spectrum always lies between the Neumann and the Dirichlet spectra \cite{bucur}:
\begin{equation} \label{comparison}
    \lambda_0^j < \lambda_\varepsilon^j < \lambda_{\infty}^j, \quad \varepsilon > 0,\quad  j \geq 0.
\end{equation}
Here, $\lambda_\infty^j$ denotes the Dirichlet eigenvalue with index $j$, starting from $j=0$.

Appendix~\ref{sec:proofs} proves the following Robin eigenvalue expansions. For the ground state,
\begin{equation}
\label{asym_lambda-fixed}
 \lambda_{\varepsilon}^{0} = \frac{\varepsilon | \partial \Omega|}{|\Omega|}  - \frac{\varepsilon^2}{|\Omega|} \sum_{j=1}^{+\infty} \frac{(\int_{\partial \Omega} u_0^j \, d\sigma)^2}{\lambda_0^j} + \mathcal{O}(\varepsilon^3) \quad \text{as } \varepsilon \to 0  
\end{equation}
and, for fixed $l\geq 1$ such that $\lambda_0^l$ is simple, 
\begin{equation} \label{asymp:eigenvaluel}
 \lambda_{\varepsilon}^{l} = \lambda^l_0 + \varepsilon \int_{\partial \Omega} (u_0^l)^2 \, d\sigma + \varepsilon^2 \sum_{j\neq l} \frac{(\int_{\partial \Omega} u_0^l  u_0^j \, d\sigma)^2}{\lambda^l_0-\lambda^j_0 } + \mathcal{O}(\varepsilon^3).
\end{equation}
If $\lambda^l_0$ is of multiplicity $m$, then choose an orthonormal basis $u_0^{l,1}, \dots, u_0^{l,m}$ of its eigenspace and define the matrix $C=(C_{ij})_{i,j=1}^m$ by
$$
C_{ij}:= \int_{\partial \Omega} u_0^{l,i} u_0^{l,j}\, d\sigma.
$$
If $\mu_1,\dots,\mu_m$ are the eigenvalues of $C$, then we have
\begin{equation} \label{asymp:eigenvaluem}
\lambda_{\varepsilon}^{l,i} = \lambda^l_0 + \varepsilon \mu_i + \mathcal{O}(\varepsilon^2), 
\end{equation}
for the $i$-th branch issuing from $\lambda^l_0$. 

For a simple Neumann eigenvalue \(\lambda_0^l\), let \(u_\varepsilon^l\) denote the real \(L^2(\Omega)\)-normalized eigenfunction of the corresponding Robin branch $\lambda_\varepsilon^l$, with the sign chosen so that $(u_0^l,u_\varepsilon^l)_{L^2(\Omega)}>0$.
The following expansion holds in $H^1(\Omega)$:
\begin{equation} \label{asymp:eigenfunction}
    u_\varepsilon^l= u_0^l + \varepsilon \sum_{j\neq l} \frac{\int_{\partial \Omega} u_0^l  u_0^j \, d\sigma}{\lambda^l_0-\lambda^j_0 }
u_0^j + \mathcal{O}(\varepsilon^2).
\end{equation}
In particular, for $l=0$, we have
\begin{equation} \label{asymp:eigenfunction0}
    u_\varepsilon^0= \frac{1}{\sqrt{|\Omega|}} - \frac{1}{\sqrt{|\Omega|}} \varepsilon \sum_{j=1}^{+\infty} \frac{\int_{\partial \Omega} u_0^j \, d\sigma}{\lambda^j_0 }
u_0^j + \mathcal{O}(\varepsilon^2).
\end{equation}

\subsection{Uniform high frequency asymptotics}
Note that the remainder constants in the fixed-index expansions above can depend on $l$.  A quantitative uniform Taylor estimate requires spectral separation. We refer the reader to Proposition~\ref{prop:gap-remainder-A} which makes this dependence explicit.

Here, we first obtain an absolute Robin--Neumann comparison without a
gap assumption. Throughout this subsection, we assume that \(\partial\Omega\) is of class \(C^\infty\). 

For all $v\in H^1(\Omega)$, we have
\begin{equation} \label{var:form}
    \int_\Omega \nabla u_\varepsilon^j \cdot \nabla v + \varepsilon \int_{\partial \Omega} u_\varepsilon^j v = \lambda_\varepsilon^j \int_{\Omega} u_\varepsilon^j v.
\end{equation}

Note that the differentiation below is performed along locally analytic normalized eigenbranches. For ordered eigenvalues that are locally absolutely continuous, the derivative identity is interpreted almost everywhere, with a corresponding choice of normalized eigenfunctions.

Taking the derivative of \eqref{var:form} with respect to $\varepsilon$ and using the fact that $\int_{\Omega} u_\varepsilon^j \frac{\mathrm{d} u_\varepsilon^j}{\mathrm{d}\varepsilon} =0$, it follows that
\begin{equation} \label{der:lambda}
    \frac{\mathrm{d} \lambda_\varepsilon^j}{\mathrm{d} \mathrm{\varepsilon}} = \int_{\partial \Omega} (u_\varepsilon^j)^2\, d\sigma. 
\end{equation}
Therefore, 
\begin{equation} \label{form:lambda}
    \lambda_\varepsilon^j = \lambda_0^j + \int_0^\varepsilon\int_{\partial \Omega} (u_s^j)^2\, d\sigma\, ds. 
\end{equation}

Let 
$$
d_j(\varepsilon):= \lambda_\varepsilon^j - \lambda_0^j, \quad j \geq 0. 
$$
From \eqref{form:lambda}, the trace inequality
$$
\int_{\partial \Omega} (u_s^j)^2\, d\sigma\, \leq C_\Omega \|u_s^j\|_{L^2(\Omega)} \, \| u_s^j\|_{H^1(\Omega)}, 
$$
and 
$$
\| u_s^j\|_{H^1(\Omega)} \leq C_\Omega (1+\lambda_s^j)^{1/2} \leq C_\Omega (1+ \lambda_\infty^j)^{1/2}   \quad (\text{by using \eqref{comparison}}),
$$
it follows that 
$$
0\leq d_j(\varepsilon) \leq C_\Omega \varepsilon (1+ \lambda_\infty^j)^{1/2},
$$
for some constant $C_\Omega$ that depends only on $\Omega$. 

By Weyl's law, $\lambda_\infty^j\asymp j^{2/d}$ as $j\to\infty$.
Consequently, we obtain
\[
 \sup_{0\leq j\leq J(\varepsilon)}
       |\lambda_\varepsilon^j-\lambda_0^j|\longrightarrow0
 \quad\hbox{if}\quad
 \varepsilon(1+J(\varepsilon))^{1/d}\longrightarrow0,
\]
that is, for $J(\varepsilon)=o(\varepsilon^{-d})$. 

Next, we turn to the asymptotic expansion \eqref{asymp:eigenvaluel}. We assume that $\lambda_0^l$ is simple and let 
\begin{equation} \label{def:gammal}
    \gamma^l:= \inf_{j\neq l} |\lambda^j_0 - \lambda_0^l| >0
\end{equation}
and
\begin{equation} \label{def:Bl} 
B_l:= \sum_{j\neq l} \frac{a_{lj}^2}{\lambda^l_0-\lambda^j_0 }, \quad a_{lj}:= \int_{\partial \Omega} u_0^l  u_0^j \, d\sigma.
\end{equation}

To prove the absolute convergence of the series defining $B_l$, let $w^l\in H^1(\Omega)$ solve
\begin{equation}\label{prob:aux}
\left\{ \begin{array}{l}
(-\Delta+1)w^l=0\ \text{in }\Omega,\\
\nm \ds 
\frac{\partial w^l}{\partial \nu} =u^l_0\ \text{on }\partial\Omega.
\end{array}
\right.
\end{equation}
Testing the weak equation of \eqref{prob:aux} with $u^j_0$ gives
\[
a_{lj}=(1+\lambda^j_0)(w^l,u_0^j)_{L^2(\Omega)}.
\]
Therefore, we obtain
\begin{equation}
\sum_j\frac{|a_{lj}|^2}{1+\lambda_0^j}
=\sum_j(1+\lambda^j_0)|(w^l,u_0^j)_{L^2(\Omega)}|^2<\infty,
\end{equation}
since $w^l\in H^1(\Omega)$ and $(u_0^j)_j$ is an orthonormal basis of $L^2(\Omega)$.

For fixed $l$, Weyl's law implies $|\lambda_0^j-\lambda_0^l|\asymp 1+\lambda_0^j$ for large $j$. Hence, we have
\begin{equation}
\sum_{j\neq l}\frac{|a_{lj}|^2}{|\lambda_0^j-\lambda_0^l|}< + \infty.
\end{equation}
Thus,  the series in \eqref{def:Bl} is absolutely convergent. 

Moreover, with $\gamma^l$ as in \eqref{def:gammal},
\[
 |B_l|\leq
 \left(\sup_{j\ne l}\frac{1+\lambda_0^j}{|\lambda_0^l-\lambda_0^j|}\right)
       \sum_{j\ne l}\frac{|a_{lj}|^2}{1+\lambda_0^j}
 \leq\left(1+\frac{1+\lambda_0^l}{\gamma^l}\right)
       \sum_{j\ne l}\frac{|a_{lj}|^2}{1+\lambda_0^j}.
\]
The lifting estimate gives
\[
 \sum_j\frac{|a_{lj}|^2}{1+\lambda_0^j}
 \leq C_\Omega\|u_0^l|_{\partial\Omega}\|_{H^{-1/2}(\partial\Omega)}^2
 \leq C_\Omega\|u_0^l\|_{L^2(\partial\Omega)}^2.
\]
For $C^\infty$-domains the sharp boundary restriction estimate is
\begin{equation}\label{tataru}
 \|u_0^l\|_{L^2(\partial\Omega)}^2
                    \leq C_\Omega(1+\lambda_0^l)^{1/3};
\end{equation}
see \cite{tataru}. Thus,
\[
 |B_l|\leq C_\Omega(1+\lambda_0^l)^{1/3}
                     \left(1+\frac{1+\lambda_0^l}{\gamma^l}\right).
\]
A bound on $B_l$ alone does not control the higher Taylor coefficients.
The contour argument in Proposition~\ref{prop:gap-remainder-A} supplies the
additional estimate
\[
 \left|\lambda_\varepsilon^l-\lambda_0^l
       -\varepsilon a_{ll}-\varepsilon^2B_l\right|
 \leq C_\Omega a_{ll}\varepsilon^3
                    \left(1+\frac{1+\lambda_0^l}{\gamma^l}\right)^2
\]
when $\varepsilon(1+(1+\lambda_0^l)/\gamma^l)$ is sufficiently small.
In particular, the Taylor remainder is $o(\varepsilon a_{ll})$ along
simple-mode sequences satisfying
\begin{equation}\label{high-frequency}
 \varepsilon\left(1+\frac{1+\lambda_0^l}{\gamma^l}\right)\longrightarrow0.
\end{equation}
At high frequency this is equivalent, up to domain-dependent constants,
to $\varepsilon l^{2/d}/\gamma^l\to0$.  This condition concerns a fixed
smooth domain and does not assert uniformity in the prefractal depth.

\subsection{Spectral decomposition of the Green function}

We have the following spectral decomposition of the Green function $G_\varepsilon$:
\begin{equation} \label{decomp:green}
G_\varepsilon(x,y)= - \sum_{j=0}^{+\infty} \frac{u_\varepsilon^j(x) u_\varepsilon^j(y)}{\lambda_\varepsilon^j} \quad \text{ in } L^2(\Omega\times \Omega). 
\end{equation}
The positive-mode series is uniformly bounded in $L^2(\Omega\times\Omega)$, since \eqref{comparison} and Weyl's law give $\sum_{j\geq1}(\lambda_\varepsilon^j)^{-2}\leq\sum_{j\geq1}(\lambda_0^j)^{-2}<\infty$ for $d=2,3$.
Using \eqref{asym_lambda-fixed} and \eqref{asymp:eigenfunction0}, we obtain in $L^2(\Omega\times\Omega)$
\begin{equation}
    \frac{u_\varepsilon^0(x) u_\varepsilon^0(y)}{\lambda_\varepsilon^0} = \frac{1}{\varepsilon |\partial \Omega|} + \mathcal{O}(1),
\end{equation}
which shows that the singular part of the Green function $G_\varepsilon$ is given by
$$
G_\varepsilon^{\mathrm{sing}}(x,y) = - \frac{1}{\varepsilon |\partial \Omega|}.
$$
Therefore, the singular part of the resolvent is rank one and is generated by the perturbation of the Neumann zero mode. In the next section, we seek estimates uniform in the prefractal depth.

\section{Finite-depth prefractal boundary}\label{sec:prefractal}

Throughout this section, all geometric and spectral quantities are non-dimensionalized using a fixed reference length; hence the arguments of all logarithms below are dimensionless.

For each prescribed boundary approximation scale $\delta>0$, let $\Omega_\delta$ be a bounded connected Lipschitz prefractal domain and set $\Gamma_\delta:=\partial\Omega_\delta$. Here, $\delta$ represents the characteristic diameter of the boundary pieces retained in the approximation; for example, $\delta=3^{-m}$ for the level-$m$ Koch snowflake approximation. Write
\[
\sigma_\delta:=\mathbf{1}_{\Gamma_{\delta}}\mathcal H^{d-1}, \qquad |\Gamma_\delta|:=\sigma_\delta(\Gamma_\delta),
\]
where $\mathcal H^{d-1}$ is the $(d-1)$-dimensional Hausdorff measure. Let $T_\delta:H^1(\Omega_\delta)\to L^2(\Gamma_\delta,\sigma_\delta)$ be the trace operator. For $\varepsilon>0$, define the bilinear form corresponding to the Robin boundary condition
\begin{equation}\label{prefractal-robin-form}
a_{\varepsilon,\delta}(u,v) :=\int_{\Omega_\delta}\nabla u\cdot\nabla v\,dx +\varepsilon\int_{\Gamma_\delta}T_\delta u\,T_\delta v\,d\sigma_\delta, \qquad u,v\in H^1(\Omega_\delta).
\end{equation}
Let $A_{\varepsilon,\delta}$ be its associated positive self-adjoint operator.  This operator has a compact resolvent, and its ground-state eigenvalue is well-defined variationally by
\begin{equation}\label{prefractal-first-eigenvalue}
\lambda^0_{\varepsilon,\delta} =\min_{u\in H^1(\Omega_\delta)\setminus\{0\}} \frac{a_{\varepsilon,\delta}(u,u)} {\|u\|_{L^2(\Omega_\delta)}^2}>0.
\end{equation}
We list the remaining eigenvalues, repeated according to multiplicity, as
\[
0<\lambda^0_{\varepsilon,\delta} <\lambda^1_{\varepsilon,\delta} \leq\lambda^2_{\varepsilon,\delta}\leq\cdots\uparrow+\infty.
\]
We denote the integral kernel of $-A_{\varepsilon,\delta}^{-1}$ by $G_{\varepsilon,\delta}$, consistently with the sign convention in \eqref{def:green}.

To separate the total boundary mass from its distribution, we normalize the boundary measure:
\begin{equation}\label{normalized-prefractal-measure}
\widehat\sigma_\delta:=|\Gamma_\delta|^{-1}\sigma_\delta, \qquad \varepsilon\sigma_\delta =(\varepsilon|\Gamma_\delta|)\widehat\sigma_\delta.
\end{equation}
Here, $\widehat\sigma_\delta$ is a probability measure.

\begin{definition}
For coupled sequences $\varepsilon_m>0$ with $\varepsilon_m\to0$ and $\delta_m\downarrow0$, recall that the total Robin boundary mass $\tau_m$ is defined by
\[
\tau_m:=\varepsilon_m|\Gamma_{\delta_m}|,\qquad \varepsilon_m\sigma_{\delta_m}=\tau_m\widehat\sigma_{\delta_m}.
\]
We call the coupling \emph{subcritical}, \emph{critical}, or \emph{supercritical} according as
\[
\tau_m\longrightarrow0, \qquad \tau_m\longrightarrow\tau_*\in(0,\infty), \qquad\text{or}\qquad \tau_m\longrightarrow+\infty,
\]
respectively.  

When a coupling is given as a function rather than a sequence, we write it as $\delta=\delta(\varepsilon)$.
\end{definition}
\subsection{Uniform estimates and bounded-mass recovery}

Fix a bounded connected target domain $\Omega$ and write $\Gamma:=\partial\Omega$. Choose approximation scales $1>\delta_1>\delta_2>\cdots\downarrow0$ satisfying
\begin{equation}\label{scale-sampling-condition}
\frac{\log(1/\delta_{m+1})}{\log(1/\delta_m)}\longrightarrow1.
\end{equation}
For each $m$, $\delta_m$ is the diameter scale of the boundary pieces retained in the $m$th approximation, and we set
\[
\Omega_m:=\Omega_{\delta_m},\quad \Gamma_m:=\Gamma_{\delta_m},\quad \widehat\sigma_m:=\widehat\sigma_{\delta_m},\quad T_m:=T_{\delta_m},\quad \Lambda_m(\varepsilon):=\lambda^0_{\varepsilon,\delta_m}.
\]
Condition~\eqref{scale-sampling-condition} ensures that the covering exponents computed along the approximation scales $\delta_m$ recover the lower and upper Minkowski dimensions defined using all sufficiently small scales.
It holds, for example, when $\delta_m=r^m$ with $0<r<1$. Next, we state the definition of Minkowski dimension.

\begin{definition}[Minkowski dimensions and nondegeneracy]\label{def:minkowski-dimensions} 
Let $E\subset\mathbb R^d$ be nonempty and compact. For $r>0$, define its covering number by 
\begin{equation}\label{covering-number} 
\mathcal N(E,r):=\min\left\{N\in\mathbb N:E\subset\bigcup_{j=1}^{N}B(x_j,r)\text{ for some }x_1,\ldots,x_N\in\mathbb R^d\right\}.
\end{equation} 
The lower and upper Minkowski dimensions are
\begin{equation}\label{minkowski-dimensions} 
D_-(E):=\liminf_{r\downarrow0}\frac{\log\mathcal N(E,r)}{\log(1/r)},\qquad D_+(E):=\limsup_{r\downarrow0}\frac{\log\mathcal N(E,r)}{\log(1/r)}.
\end{equation} 
If these agree, their common value is denoted by $\dim_M E$. We call $E$ Minkowski nondegenerate of dimension $D$ if 
\begin{equation}\label{regular-covering-law} 
\mathcal N(E,r)\asymp r^{-D}\qquad\text{as }r\downarrow0. 
\end{equation} 
Here and below, $X\asymp Y$ means that their ratio is bounded above and away from zero by constants independent of the scale. This condition implies $\dim_M E=D$ and is stronger than the existence of the Minkowski dimension alone.
\end{definition}

\begin{remark}[Equivalent tubular-neighborhood formulation] 
Let $E^r:=\{x\in\mathbb R^d:\operatorname{dist}(x,E)\leq r\}$ and $V_E(r):=\mathcal L^d(E^r)$. 
Then we have
\begin{equation}\label{tube-covering-comparison} V_E(r)\asymp r^d\mathcal N(E,r).
\end{equation}
This follows from covering and packing estimates. An $r$-ball cover of $E$ becomes a $2r$-ball cover of $E^r$, giving the upper bound for $V_E(r)$. Conversely, a maximal $r$-separated subset of $E$ has at least $\mathcal N(E,r)$ points, whose disjoint $r/2$-balls lie in $E^r$, giving the lower bound.

Consequently, \eqref{minkowski-dimensions} agrees with the classical Euclidean specialization of the tubular-neighborhood definitions in \cite[Section~2.1]{Chow2025Fractality}. Moreover, \eqref{regular-covering-law} is equivalent to \[ 0<\liminf_{r\downarrow0}r^{D-d}V_E(r)\leq\limsup_{r\downarrow0}r^{D-d}V_E(r)<\infty. \] Thus the covering formulation gives the usual notion of Minkowski nondegeneracy. Minkowski measurability additionally requires $r^{D-d}V_E(r)$ to converge to a positive finite limit.
\end{remark}

We place the bounded connected Lipschitz domains $\Omega_m$ and the target domain $\Omega$ in a fixed bounded Lipschitz ambient domain $\mathcal O$, with $\overline\Omega\cup\bigcup_m\overline{\Omega_m}\subset\mathcal O$, and assume the volumes of the approximating domains are bounded above and away from zero
\begin{equation}\label{uniform-volume-bounds}
0<\inf_m|\Omega_m|\leq\sup_m|\Omega_m|<\infty.
\end{equation}
The following geometric conditions are sufficient for the dimension recovery estimates below, with constants independent of $m$.
\begin{enumerate}[label=\textnormal{(R\arabic*)}]
\item \emph{Uniform domains.} There exists $\eta\in(0,1)$, independent of $m$, such that every $\Omega_m$ is an $(\eta,\infty)$-domain: any distinct $x,y\in\Omega_m$ can be joined by a rectifiable curve $\gamma\subset\Omega_m$ satisfying \[ \ell(\gamma)\leq\eta^{-1}|x-y|,\qquad \operatorname{dist}(z,\partial\Omega_m)\geq\eta\frac{|x-z|\,|y-z|}{|x-y|},\qquad z\in\gamma. \]

\item \emph{Uniform upper growth bound for the normalized boundary measures.} For some $d-2<s\leq d-1$,
\begin{equation}\label{normalized-boundary-growth}
\widehat\sigma_m(B(x,r))\leq C_\sigma r^s,\qquad x\in\mathbb R^d,\quad 0<r\leq1.
\end{equation}

\item \emph{Boundary resolution.} The boundary measure at scale $\delta_m$ satisfies
\begin{equation}\label{boundary-covering-comparison}
|\Gamma_m|\asymp\delta_m^{d-1}\mathcal N(\Gamma,\delta_m).
\end{equation}
\end{enumerate}

Condition \textnormal{(R1)} ensures that the interior domain is uniformly connected, while \textnormal{(R2)} controls normalized boundary mass in small balls with radius $r$. The exponent \(s\) controls the concentration of the finite-stage boundary measures and need not equal the Minkowski dimension of the limiting boundary. For the Koch approximations, the uniform bound holds with \(s=1\), although the limiting dimension is \(D=\log4/\log3>1\). Together, these assumptions yield the normalized trace estimate used below. 

A sufficient way to verify \textnormal{(R3)} is to decompose $\Gamma_m$ into $M_m$ boundary pieces $F_{m,k}$ with pairwise overlaps of zero $\mathcal H^{d-1}$-measure, satisfying $\operatorname{diam}F_{m,k}\asymp\delta_m$, $\mathcal H^{d-1}(F_{m,k})\asymp\delta_m^{d-1}$, and $M_m\asymp\mathcal N(\Gamma,\delta_m)$, uniformly in $m$ and $k$. Summing their measures gives $|\Gamma_m|\asymp M_m\delta_m^{d-1}$, and hence \eqref{boundary-covering-comparison}. If $\Gamma$ is also Minkowski nondegenerate of dimension $D\in(d-1,d)$, then \eqref{regular-covering-law} yields \begin{equation}\label{regular-prefractal-perimeter} |\Gamma_m|\asymp\delta_m^{d-1-D}. \end{equation} Proposition~\ref{prop:koch-assumptions} verifies these comparisons for the Koch snowflake by a covering and packing argument.

We next derive the uniform estimates used in the proofs. By \textnormal{(R1)}, \eqref{uniform-volume-bounds}, and the Sobolev extension theorem \cite{Jones1981SobolevExtension}, there are linear operators $E_m:H^1(\Omega_m)\to H^1(\mathcal O)$ satisfying
\begin{equation}\label{operator-extension-bounds}
(E_mv)|_{\Omega_m}=v,\qquad \|E_mv\|_{H^1(\mathcal O)}\leq C_E\|v\|_{H^1(\Omega_m)}.
\end{equation}
The measure trace theorem and \eqref{normalized-boundary-growth} then give
\[
\|T_mv\|_{L^2(\Gamma_m,\widehat\sigma_m)}^2\leq C\left(\|v\|_{L^2(\Omega_m)}^2+\|\nabla v\|_{L^2(\Omega_m)}^2\right),\qquad v\in H^1(\Omega_m).
\]
The common uniform-domain geometry and bounded ambient domain also give the mean-zero Poincar\'e estimate
\[
\|w\|_{L^2(\Omega_m)}^2\leq C\|\nabla w\|_{L^2(\Omega_m)}^2,\qquad w\in H^1(\Omega_m),\quad \int_{\Omega_m}w\,dx=0.
\]
Combining these estimates yields
\begin{equation}\label{uniform-prefractal-estimates}
\|w\|_{L^2(\Omega_m)}^2+\|T_mw\|_{L^2(\Gamma_m,\widehat\sigma_m)}^2\leq C\|\nabla w\|_{L^2(\Omega_m)}^2,\qquad w\in H^1(\Omega_m),\quad \int_{\Omega_m}w\,dx=0;
\end{equation}
see \cite{Biegert2009Traces,HinzRozanovaTeplyaev2023}.

\begin{lemma}[Uniform ground-state eigenvalue bounds]
\label{lem:uniform-prefractal-ground-state}
Assume \textnormal{(R1)}--\textnormal{(R2)}. There is $c_0>0$, independent of $m$, such that, for every $m$, $\varepsilon>0$, 
\begin{equation}\label{bounded-mass-comparison}
c_0\min\{\varepsilon|\Gamma_m|,1\}\leq\Lambda_m(\varepsilon)\leq\frac{\varepsilon|\Gamma_m|}{|\Omega_m|}.
\end{equation}

Moreover, there are $C,\tau_0>0$, independent of $m$, such that
\begin{equation}\label{uniform-prefractal-slope}
0\leq\frac{\varepsilon|\Gamma_m|}{|\Omega_m|}-\Lambda_m(\varepsilon) \leq C(\varepsilon|\Gamma_m|)^2, \qquad 0<\varepsilon|\Gamma_m|\leq\tau_0.
\end{equation}
Thus, uniformly in $m$,
\begin{equation}\label{uniform-prefractal-expansion}
\Lambda_m(\varepsilon) =\frac{\varepsilon|\Gamma_m|}{|\Omega_m|} \left(1+\mathcal O(\varepsilon|\Gamma_m|)\right) \qquad\text{as }\varepsilon|\Gamma_m|\to0.
\end{equation}
\end{lemma}

\begin{proof}
Let $\tau := \varepsilon|\Gamma_m|$ and $u=\bar u+w\in H^1(\Omega_m)$ be any $L^2(\Omega_m)$-normalized function, where $\bar u:=|\Omega_m|^{-1}\int_{\Omega_m}u\,dx$ and $\int_{\Omega_m}w\,dx=0$. Since $\widehat\sigma_m$ is a probability measure and $\bar u = T_m u - T_m w$
\[
|\bar u|\leq \|T_mu\|_{L^2(\Gamma_m,\widehat\sigma_m)} +\|T_mw\|_{L^2(\Gamma_m,\widehat\sigma_m)},
\]
and hence, by \eqref{uniform-prefractal-estimates} and \eqref{uniform-volume-bounds},
\begin{align*}
1=& \|u\|_{L^2(\Omega_m)}^2 = |\Omega_m|\bar u^2+\|w\|_{L^2(\Omega_m)}^2\leq C(\|T_mu\|_{L^2(\Gamma_m,\widehat\sigma_m)}^2+\|\nabla w\|_{L^2(\Omega_m)}^2)\\
\leq& C\max\{1,\tau^{-1}\}(\tau \|T_mu\|_{L^2(\Gamma_m,\widehat\sigma_m)}^2+\|\nabla w\|_{L^2(\Omega_m)}^2) = C\max\{1,\tau^{-1}\} a_{\varepsilon,\delta_m}(u,u).
\end{align*}
Taking the infimum proves the lower bound in \eqref{bounded-mass-comparison}, after decreasing $c_0$ if necessary; the constant test function proves the upper bound.

Next, take $u$ to be the positive normalized ground state. Testing its weak eigenvalue equation with $w$ and using $T_mu=\bar u+T_mw$ gives the exact identity
\[
\|\nabla w\|_{L^2(\Omega_m)}^2+\tau\|T_mw\|_{L^2(\Gamma_m,\widehat\sigma_m)}^2 +\tau\bar u\int_{\Gamma_m}T_mw\,d\widehat\sigma_m =\Lambda_m(\varepsilon)\|w\|_{L^2(\Omega_m)}^2.
\]
Using \eqref{uniform-prefractal-estimates}, the fact that $\widehat\sigma_m$ is a probability measure, and the upper bound in \eqref{bounded-mass-comparison}, we obtain
\[
g\leq C\tau |\bar u|\sqrt g+\frac{C\tau}{|\Omega_m|}g,
\]
where $g:=\|\nabla w\|_{L^2(\Omega_m)}^2$. 
Normalization and \eqref{uniform-volume-bounds} give $|\bar u|\leq|\Omega_m|^{-1/2}\leq C$. For sufficiently small $\tau$, the last term is absorbed into the left-hand side. If $g>0$, division by $\sqrt g$ then gives $\sqrt g\leq C\tau|\bar u|$; the case $g=0$ is immediate. Thus, $g=\mathcal O(\tau^2)$ uniformly in $m$, and \eqref{uniform-prefractal-estimates} also gives
\[
\|w\|_{L^2(\Omega_m)}+\|T_mw\|_{L^2(\Gamma_m,\widehat\sigma_m)}=\mathcal O(\tau).
\]
Normalization now yields
\[
\bar u^2=|\Omega_m|^{-1}+\mathcal O(\tau^2),
\]
and expanding the normalized boundary energy gives
\[
 \|\bar u+T_mw\|_{L^2(\Gamma_m,\widehat\sigma_m)}^2 =\bar u^2+2\bar u\int_{\Gamma_m}T_mw\,d\widehat\sigma_m+\|T_mw\|_{L^2(\Gamma_m,\widehat\sigma_m)}^2 =|\Omega_m|^{-1}+\mathcal O(\tau).
\]
Therefore $\Lambda_m(\varepsilon)=g+\tau \|T_mu\|_{L^2(\Gamma_m,\widehat\sigma_m)}^2=\tau/|\Omega_m|+\mathcal O(\tau^2)$. Together with the constant-function upper bound $\Lambda_m(\varepsilon)\leq\tau/|\Omega_m|$, this proves \eqref{uniform-prefractal-slope}.
\end{proof}

\begin{cor}[Subcritical Green's singularity]
\label{cor:subcritical-green-singularity}
Assume \textnormal{(R1)}--\textnormal{(R2)}. Let $P_{0,m}$ be the orthogonal projection of $L^2(\Omega_m)$ onto the constant functions, defined by
\[
P_{0,m}f:=\frac1{|\Omega_m|}\left(\int_{\Omega_m}f\,dx\right)\mathbf 1_{\Omega_m}.
\]
There are $C,\tau_0>0$, independent of $m$, such that, whenever $0<\tau:=\varepsilon|\Gamma_m|\leq\tau_0$,
\begin{equation}\label{subcritical-rank-one-resolvent}
\left\|A_{\varepsilon,\delta_m}^{-1}-\frac{|\Omega_m|}{\tau}P_{0,m}\right\|_{\mathcal L(L^2(\Omega_m))}\leq C,
\end{equation}
where $\mathcal L(L^2(\Omega_m))$ is the space of bounded linear operators from $L^2(\Omega_m)$ into $L^2(\Omega_m)$.
Consequently, along every subcritical sequence $\tau_m\to0$, the signed Green operator has leading term $-(|\Omega_m|/\tau_m)P_{0,m}$. Since the kernel of $P_{0,m}$ is $|\Omega_m|^{-1}$, its singular rank-one kernel is equal to
\[
-\frac1{\tau_m}=-\frac1{\varepsilon_m|\Gamma_m|}.
\]
\end{cor}

\begin{proof}
For fixed $m$ and $\tau$, let $u$ be the positive normalized ground state and $P$ the orthogonal projection onto its span. The proof of Lemma~\ref{lem:uniform-prefractal-ground-state} gives
\[
\Lambda_m(\varepsilon)=\frac{\tau}{|\Omega_m|}+\mathcal O(\tau^2), \qquad \|u-|\Omega_m|^{-1/2}\mathbf 1_{\Omega_m}\|_{L^2(\Omega_m)}=\mathcal O(\tau),
\]
uniformly in $m$. Hence,
\[
\Lambda_m(\varepsilon)^{-1}=\frac{|\Omega_m|}{\tau}+\mathcal O(1), \qquad \|P-P_{0,m}\|_{\mathcal L(L^2(\Omega_m))}=\mathcal O(\tau).
\]
The mean-zero Poincar\'e estimate in \eqref{uniform-prefractal-estimates} gives a uniform lower bound for the first positive Neumann eigenvalue. The min--max principle then gives
\[
\lambda^1_{\varepsilon,\delta_m}\geq\lambda^1_{0,\delta_m}\geq C^{-1}.
\]
The spectral theorem, therefore, gives
\[
\left\|A_{\varepsilon,\delta_m}^{-1}-\Lambda_m(\varepsilon)^{-1}P\right\|_{\mathcal L(L^2(\Omega_m))}\leq C.
\]
Combining the preceding estimates explicitly,
\[
\left\|A_{\varepsilon,\delta_m}^{-1}-\frac{|\Omega_m|}{\tau}P_{0,m}\right\| \leq C+\left|\Lambda_m(\varepsilon)^{-1}-\frac{|\Omega_m|}{\tau}\right|+\frac{|\Omega_m|}{\tau}\|P-P_{0,m}\|\leq C,
\]
which proves \eqref{subcritical-rank-one-resolvent}.  
\end{proof}

\begin{thm}[Spectral recovery of the lower and upper Minkowski dimensions]
\label{thm:minkowski-recovery}
Assume \textnormal{(R1)}--\textnormal{(R3)} and \eqref{scale-sampling-condition}. Let $\varepsilon_m>0$ satisfy $\varepsilon_m\to0$ and suppose $\sup_m\tau_m<\infty$. Then
\begin{equation}\label{bounded-mass-perimeter-equivalence}
\frac{\Lambda_m(\varepsilon_m)}{\varepsilon_m}\asymp|\Gamma_m|.
\end{equation}
Consequently,
\begin{equation}\label{lower-minkowski-recovery}
D_-(\Gamma)=d-1+\liminf_{m\to\infty}\frac{\log\!\left(\Lambda_m(\varepsilon_m)/\varepsilon_m\right)}{\log(1/\delta_m)},
\end{equation}
\begin{equation}\label{upper-minkowski-recovery}
D_+(\Gamma)=d-1+\limsup_{m\to\infty}\frac{\log\!\left(\Lambda_m(\varepsilon_m)/\varepsilon_m\right)}{\log(1/\delta_m)}.
\end{equation}
If $\dim_M\Gamma$ exists, the liminf and limsup are equal and may be replaced by the ordinary limit. In the subcritical regime $\tau_m\to0$, we have the stronger form
\begin{equation}\label{subcritical-exact-perimeter}
\frac{\Lambda_m(\varepsilon_m)}{\varepsilon_m}=\frac{|\Gamma_m|}{|\Omega_m|}\bigl(1+o(1)\bigr)
\end{equation}

\end{thm}

\begin{proof}
Let $B:=\sup_m\tau_m$. Dividing \eqref{bounded-mass-comparison} by $\varepsilon_m=\tau_m/|\Gamma_m|$ and using \eqref{uniform-volume-bounds} gives
\[
\frac{c_0}{\max\{1,B\}}|\Gamma_m|\leq\frac{\Lambda_m(\varepsilon_m)}{\varepsilon_m}\leq C|\Gamma_m|.
\]
Together with \eqref{boundary-covering-comparison}, this yields
\begin{equation}\label{eq:spectral-covering-log}
\log\!\left(\frac{\Lambda_m(\varepsilon_m)}{\varepsilon_m}\right)=\log\mathcal N(\Gamma,\delta_m)-(d-1)\log(1/\delta_m)+\mathcal O(1).
\end{equation}
If $\delta_{m+1}<r\leq\delta_m$, monotonicity of the covering number gives
\[
\frac{\log\mathcal N(\Gamma,\delta_m)}{\log(1/\delta_{m+1})}\leq\frac{\log\mathcal N(\Gamma,r)}{\log(1/r)}\leq\frac{\log\mathcal N(\Gamma,\delta_{m+1})}{\log(1/\delta_m)}.
\]
The ratios of the two neighboring logarithmic scales tend to one by \eqref{scale-sampling-condition}. Moreover, compactness of $\Gamma$ gives $\mathcal N(\Gamma,r)\leq Cr^{-d}$ for small $r$, so both bounding quotients are bounded. Taking lower and upper limits and using \eqref{scale-sampling-condition} identifies the exponents computed along the approximation scales with \eqref{minkowski-dimensions}.
The term $\mathcal O(1)/\log(1/\delta_m)$ from \eqref{eq:spectral-covering-log} tends to zero, proving \eqref{lower-minkowski-recovery}--\eqref{upper-minkowski-recovery}. Finally, \eqref{subcritical-exact-perimeter} follows from \eqref{uniform-prefractal-expansion}.
\end{proof}

We also consider reconstruction at a fixed eigenvalue level. Assume \textnormal{(R1)}--\textnormal{(R2)}, let $c_0$ be the constant in Lemma~\ref{lem:uniform-prefractal-ground-state}, and fix $\ell_*\in(0,c_0)$. The variational characterization shows that $\varepsilon\mapsto\Lambda_m(\varepsilon)$ is continuous and strictly increasing; strictness follows because the positive ground state has nonzero boundary trace. Since $\Lambda_m(0)=0$ and \eqref{bounded-mass-comparison} gives $\Lambda_m(|\Gamma_m|^{-1})\geq c_0$, there is a unique $\varepsilon_m^{\mathrm c}>0$ such that $\Lambda_m(\varepsilon_m^{\mathrm c})=\ell_*$. The same comparison yields
\begin{equation}\label{eq:fixed-level-mass}
\ell_*|\Omega_m|\leq\varepsilon_m^{\mathrm c}|\Gamma_m|\leq\frac{\ell_*}{c_0}<1.
\end{equation}
Thus $\varepsilon_m^{\mathrm c}|\Gamma_m|\asymp1$ and $\log(1/\varepsilon_m^{\mathrm c})=\log|\Gamma_m|+\mathcal O(1)$. If \textnormal{(R3)} and \eqref{scale-sampling-condition} also hold, then
\begin{equation}\label{eq:fixed-level-dimensions}
D_-(\Gamma)=d-1+\liminf_{m\to\infty}\frac{\log(1/\varepsilon_m^{\mathrm c})}{\log(1/\delta_m)},\qquad D_+(\Gamma)=d-1+\limsup_{m\to\infty}\frac{\log(1/\varepsilon_m^{\mathrm c})}{\log(1/\delta_m)}.
\end{equation}
The corresponding total Robin boundary masses are bounded above and away from zero; convergence to a critical mass is not asserted. If, moreover, $\Gamma$ is Minkowski nondegenerate of dimension $D\in(d-1,d)$, then
\[
\varepsilon_m^{\mathrm c}\asymp\delta_m^{D-(d-1)},\qquad D=d-1+\lim_{m\to\infty}\frac{\log(1/\varepsilon_m^{\mathrm c})}{\log(1/\delta_m)}.
\]

We next classify power-law couplings. Assume \textnormal{(R1)}--\textnormal{(R3)}, set $\theta_-:=D_-(\Gamma)-(d-1)$ and $\theta_+:=D_+(\Gamma)-(d-1)$, and choose $\varepsilon_m:=\delta_m^{1/a}$ with $a>0$. Under \eqref{scale-sampling-condition}, we have
\begin{equation}\label{eq:power-mass-exponents}
\liminf_{m\to\infty}\frac{\log\tau_m}{\log(1/\delta_m)}=\theta_--\frac1a,\qquad \limsup_{m\to\infty}\frac{\log\tau_m}{\log(1/\delta_m)}=\theta_+-\frac1a.
\end{equation}
With $1/0:=+\infty$, the coupling is subcritical if $a<1/\theta_+$ and supercritical if $a>1/\theta_-$. If $1/\theta_+<a<1/\theta_-$, then $\tau_m$ has subsequences tending to zero and to $+\infty$, so the total mass has no limit; the endpoint behavior is not determined by $D_-(\Gamma)$ and $D_+(\Gamma)$ alone. If $\Gamma$ is Minkowski nondegenerate of dimension $D\in(d-1,d)$, then $\theta_-=\theta_+=:\theta>0$.  Consequently, the transition window reduces to the single exponent $a_{\mathrm c}:=1/\theta$: the regimes are subcritical for $a<a_{\mathrm c}$ and supercritical for $a>a_{\mathrm c}$, while at $a=a_{\mathrm c}$ the masses remain bounded above and away from zero. This endpoint is critical in the strict sense only if $\delta_m^\theta|\Gamma_m|$ converges to a value in $(0,\infty)$.

Finally, we verify that the three assumptions hold for the standard polygonal approximations to the Koch snowflake.
\begin{prop}[The standard Koch snowflake]\label{prop:koch-assumptions}
Let $\Omega_m$ be the standard nested polygonal approximations to the planar Koch snowflake $\Omega$, starting from an equilateral triangle of side length one, and set
\[
\delta_m=3^{-m},\qquad M_m=3\cdot4^m,\qquad P_m:=|\Gamma_m|=M_m\delta_m.
\]
Then \eqref{uniform-volume-bounds}, \eqref{scale-sampling-condition}, and \textnormal{(R1)}--\textnormal{(R3)} hold with $s=1$. Moreover, $\Gamma$ is Minkowski nondegenerate of dimension $D=\log4/\log3$. Consequently, Theorem~\ref{thm:minkowski-recovery} applies to this family.
\end{prop}

\begin{proof}
Nestedness gives $0<|\Omega_1|\leq|\Omega_m|\leq|\Omega|<\infty$, while $\log(1/\delta_{m+1})/\log(1/\delta_m)=(m+1)/m\to1$, proving \eqref{uniform-volume-bounds} and \eqref{scale-sampling-condition}.

Condition \textnormal{(R1)} follows from \cite[Proof of Proposition~3.2]{Pang2008NeumannPart2}, which establishes the $(\eta,\infty)$-domain property for the standard inner Koch polygons with $\eta>0$ independent of the stage. Rescaling the initial triangle to unit side length preserves this property.

For \textnormal{(R2)}, if $0<r\leq\delta_m$, the triangular lattice of spacing $\delta_m$ implies that $B(x,r)$ meets only a bounded number of edges, each contributing length at most $2r$. Hence
\[
\widehat\sigma_m(B(x,r))\leq C\frac{r}{P_m}\leq Cr.
\]
If $\delta_m\leq r\leq1$, choose $0\leq k\leq m$ with $\delta_k\leq r<3\delta_k$. Descendants of a level-$k$ edge stay within distance $C\delta_k$ of that edge. The enlarged ball $B(x,r+C\delta_k)$ meets only a bounded number of level-$k$ lattice edges, so $B(x,r)$ meets only a bounded number of their descendant families. Each family has normalized mass $4^{m-k}\delta_m/P_m=(3\cdot4^k)^{-1}$, giving
\[
\widehat\sigma_m(B(x,r))\leq C4^{-k}=C\delta_k^D\leq Cr^D\leq Cr,
\]
with constants independent of $m$, $x$, and $r$. Thus \textnormal{(R2)} holds with $s=1$.

For \textnormal{(R3)}, each level-$m$ edge has length $\delta_m$, and its limiting descendant arc has diameter at most $C\delta_m$. These arcs cover $\Gamma$, and each can be covered by a bounded number of $\delta_m$-balls, giving $\mathcal N(\Gamma,\delta_m)\leq C_1M_m$. Conversely, the $M_m$ oriented starting vertices belong to $\Gamma$ and are distinct and $\delta_m$-separated. Each $\delta_m$-ball therefore contains only a bounded number of them, giving $\mathcal N(\Gamma,\delta_m)\geq c_1M_m$. Consequently,
\[
\mathcal N(\Gamma,\delta_m)\asymp M_m=3\delta_m^{-D}.
\]
Together with $P_m=M_m\delta_m$, this proves \textnormal{(R3)}. Since $\delta_{m+1}=\delta_m/3$, monotonicity extends the covering estimate to all sufficiently small $r>0$, proving \eqref{regular-covering-law}.
\end{proof}

\begin{figure}[H]
\centering
\includegraphics[width=0.98\textwidth]{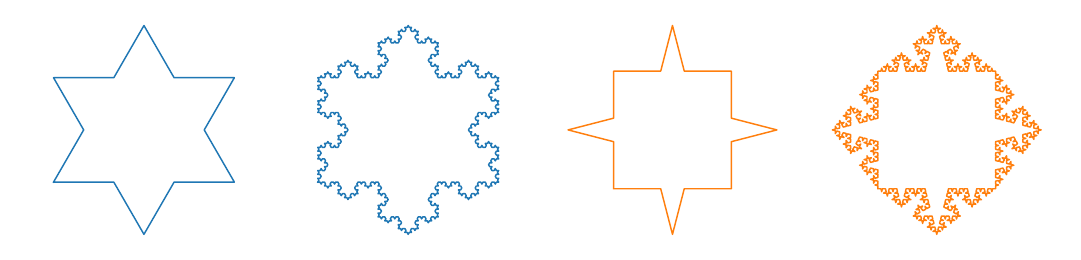}
\caption{Polygonal snowflake approximations.
From left to right: the Koch snowflake (blue) at stages $m=1$ and $m=4$, followed by the homogeneous Rohde snowflake \cite[Section~3]{Rohde2001Quasicircles} (orange) with $p=0.4$ at stages $m=1$ and $m=4$.
The Koch family has $\delta_m=3^{-m}$ and limiting boundary dimension $D=\log4/\log3\approx1.262$; the Rohde family has $\delta_m=p^m$ and limiting boundary dimension $D=\log4/\log(1/p)\approx1.513$.
The seeds at $m=0$ are a unit-side equilateral triangle and a unit-side square, respectively; $m$ counts boundary replacement steps.}
\label{fig:snowflake-geometry}
\end{figure}

Figure~\ref{fig:snowflake-geometry} illustrates these Koch approximations alongside homogeneous Rohde approximations with $p=0.4$. The latter also satisfy \textnormal{(R1)}--\textnormal{(R3)}: the uniform quasicircle bound in \cite[Lemma~3.1]{Rohde2001Quasicircles} yields \textnormal{(R1)}, while the corresponding vertex-separation estimates replace the triangular-lattice argument in Proposition~\ref{prop:koch-assumptions} and give \textnormal{(R2)}, with $s=1$, and \textnormal{(R3)}. Here, homogeneous means that every edge is replaced at every generation by the four-segment Rohde generator with contraction ratio \(p\). Starting from a unit-side square, this gives \(|\Gamma_m|=4(4p)^m\) and \(\delta_m=p^m\).

\subsection{Operator interpretation of the three regimes}\label{sec:operator-regimes}

We describe the operator limits under additional convergence hypotheses, which are not needed for the dimension recovery established above.
Throughout this subsection, the bounded connected Lipschitz domains $\Omega_m$ form a nested inner approximation of $\Omega$:
\begin{equation}\label{operator-inner-exhaustion}
\Omega_m\subset\Omega_{m+1}\subset\Omega,\qquad \bigcup_{m=1}^{\infty}\Omega_m=\Omega,
\end{equation}
with $\overline\Omega\subset\mathcal O$. The total Robin boundary mass $\tau_m=\varepsilon_m|\Gamma_m|$ does not by itself determine a limiting operator. Assume additionally that, as measures on $\overline{\mathcal O}$,
\begin{equation}\label{operator-measure-assumptions}
\widehat\sigma_m\stackrel{*}{\rightharpoonup}\widehat\sigma,\qquad \widehat\sigma(\Gamma)=1,
\end{equation}
where $\widehat\sigma$ induces a bounded trace $T:H^1(\Omega)\to L^2(\Gamma,\widehat\sigma)$ agreeing with boundary restriction for smooth functions.

Write $a_m:=a_{\varepsilon_m,\delta_m}$, and let $J_m$ and $J$ denote zero extension from $\Omega_m$ and $\Omega$ to $\mathcal O$, respectively. Regard the quadratic energies as functionals on $L^2(\mathcal O)$, with value $+\infty$ outside their zero-extended form domains. We retain the uniform extension bound \eqref{operator-extension-bounds}. Writing $\widetilde u_m:=J_mu_m$ and $\widetilde u:=Ju$, we assume Mosco convergence to the candidate limit form $a_{\mathrm{lim}}$ specified for the relevant mass regime in cases \textnormal{(i)}--\textnormal{(iii)} below. Explicitly, we require the following two conditions. First, for every sequence $u_m\in H^1(\Omega_m)$ and every $u\in L^2(\Omega)$,
\begin{equation}
\label{mosco_liminf}
\widetilde u_m\rightharpoonup\widetilde u\ \text{in }L^2(\mathcal O)\quad\Longrightarrow\quad a_{\mathrm{lim}}(u,u)\leq\liminf_{m\to\infty}a_m(u_m,u_m), 
\end{equation}
where $a_{\mathrm{lim}}(u,u):=+\infty$ if $u\notin D(a_{\mathrm{lim}})$. Conversely, for every $u\in D(a_{\mathrm{lim}})$, there exists a sequence $u_m\in H^1(\Omega_m)$ such that 
\begin{equation}
\label{mosco_limsup}
\widetilde u_m\to\widetilde u\ \text{in }L^2(\mathcal O),\qquad \limsup_{m\to\infty}a_m(u_m,u_m)\leq a_{\mathrm{lim}}(u,u). 
\end{equation}
Compactness follows from the uniform extension bound \eqref{operator-extension-bounds} and the nested exhaustion \eqref{operator-inner-exhaustion}. Indeed, if 
$\sup_m(\|u_m\|_{L^2(\Omega_m)}^2+a_m(u_m,u_m))<\infty$, then 
\[ 
\|E_mu_m\|_{H^1(\mathcal O)}^2\leq C_E^2\bigl(\|u_m\|_{L^2(\Omega_m)}^2+a_m(u_m,u_m)\bigr). \] 
Rellich compactness gives $E_mu_m\to v$ in $L^2(\mathcal O)$ along a subsequence. By \eqref{operator-inner-exhaustion}, 
\[ 
\|\widetilde u_m-\mathbf1_\Omega v\|_{L^2(\mathcal O)}\leq\|E_mu_m-v\|_{L^2(\mathcal O)}+\|(\mathbf1_{\Omega_m}-\mathbf1_\Omega)v\|_{L^2(\mathcal O)}\longrightarrow0. 
\] 
Thus, bounded norm and energy imply relative compactness of the zero extensions, so the Mosco convergence is compact.

The limit domain also admits $H^1$ extensions. For $u\in H^1(\Omega)$, a weak subsequential limit of the uniformly bounded sequence $E_m(u|_{\Omega_m})$ equals $u$ on $\Omega$ by nested exhaustion.

We briefly indicate how the assumed Mosco convergence in \eqref{mosco_liminf} and \eqref{mosco_limsup} can be verified for a specific approximation family. If $\tau_m\to\tau_*\in[0,\infty)$, uniform extensions and compactness identify limits of sequences with bounded $L^2$ norm and energy. For a subsequential limit $u$, the key boundary step is to establish 
\begin{equation}\label{eq:boundary-trace-convergence}
\int_{\Gamma_m}|T_mu_m|^2\,d\widehat\sigma_m\longrightarrow\int_\Gamma|Tu|^2\,d\widehat\sigma.
\end{equation}
If \textnormal{(R2)} holds, applying \cite[Lemma~2.1]{HinzRozanovaTeplyaev2023} to a subsequence $E_mu_m\rightharpoonup U$ in $H^1(\mathcal O)$ gives \eqref{eq:boundary-trace-convergence}, with $u=U|_\Omega$. Multiplication by $\tau_m\to\tau_*$ and weak lower semicontinuity of the interior energy then give the Mosco lower bound. For recovery, take a fixed extension $U\in H^1(\mathcal O)$ of $u$ and set $u_m:=U|_{\Omega_m}$: nested exhaustion gives convergence of the interior energies, and \eqref{eq:boundary-trace-convergence} gives convergence of the boundary energies.

If $\tau_m\to\infty$, bounded energy gives
\[
\int_{\Gamma_m}|T_mu_m|^2\,d\widehat\sigma_m\leq\frac{C}{\tau_m}\longrightarrow0.
\]
For the Mosco lower bound, the required boundary step is to establish
\[
\int_\Gamma|Tu|^2\,d\widehat\sigma\leq\liminf_m\int_{\Gamma_m}|T_mu_m|^2\,d\widehat\sigma_m=0,
\]
which forces $Tu=0$ and supplies the domain constraint in the Mosco lower bound. Recovery must also be verified. If $\ker T=H^1_0(\Omega)$, as in \eqref{operator-zero-trace-kernel}, choose $\phi_j\in C_c^\infty(\Omega)$ with $\phi_j\to u$ in $H^1(\Omega)$. Nested exhaustion permits a choice $j(m)\to\infty$ with $\operatorname{supp}\phi_{j(m)}\Subset\Omega_m$. Thus $u_m:=\phi_{j(m)}|_{\Omega_m}$ satisfies $T_mu_m=0$ and $a_m(u_m,u_m)\to\int_\Omega|\nabla u|^2\,dx$, independently of the growth rate of $\tau_m$. Without this kernel identity, recovery in the limiting zero-trace space remains to be checked.

Under these additional variational hypotheses, the three regimes have the following operator interpretation:
\begin{enumerate}[label=\textnormal{(\roman*)}]
\item If $\tau_m\to0$, the limiting form is
\[
a_{\mathrm{lim}}(u,v)=\int_\Omega\nabla u\cdot\nabla v\,dx,\qquad u,v\in H^1(\Omega),
\]
and the limit is the Neumann Laplacian.
\item If $\tau_m\to\tau_*\in(0,\infty)$, then $\varepsilon_m\sigma_{\delta_m}=\tau_m\widehat\sigma_m\stackrel{*}{\rightharpoonup}\tau_*\widehat\sigma$, and the limiting form is
\[
a_{\mathrm{lim}}(u,v)=\int_\Omega\nabla u\cdot\nabla v\,dx+\tau_*\int_\Gamma Tu\,Tv\,d\widehat\sigma,\qquad u,v\in H^1(\Omega).
\]
The associated operator is the generalized Robin Laplacian for the boundary measure $\tau_*\widehat\sigma$.
\item If $\tau_m\to+\infty$, total-mass divergence alone does not determine the limit. If compact Mosco convergence holds to
\[
a_{\mathrm{lim}}(u,v)=\int_\Omega\nabla u\cdot\nabla v\,dx,\qquad u,v\in\{w\in H^1(\Omega):Tw=0\ \widehat\sigma\text{-a.e.}\},
\]
the limit is the Laplacian associated with this zero-trace form. It is the usual Dirichlet Laplacian precisely when
\begin{equation}\label{operator-zero-trace-kernel}
\{w\in H^1(\Omega):Tw=0\ \widehat\sigma\text{-a.e.}\}=H^1_0(\Omega).
\end{equation}
\end{enumerate}
In each case where compact Mosco convergence holds to the indicated form, let $A_m$ and $A_{\mathrm{lim}}$ be the associated operators. Then, for every $\alpha>0$, the zero-extended resolvents converge in operator norm:
\begin{equation}
\label{operator-norm-resolvent-limit}
 \left\|J_m(A_m+\alpha I)^{-1}J_m^*-J(A_{\mathrm{lim}}+\alpha I)^{-1}J^*\right\|_{\mathcal L(L^2(\mathcal O))}\longrightarrow0,
 \end{equation}
and the eigenvalues converge at every fixed index, with multiplicities counted. See also \cite[Theorem~6.1]{HinzRozanovaTeplyaev2023}.

\section{Applications and numerical illustrations}\label{sec:applications}
\subsection{Scaling for self-similar boundaries}

Let $0\leq t_0 < t_1 <\dots < t_j < t_{j+1} < \dots \rightarrow +\infty $ be a time sequence. Consider self-similar boundaries and suppose that each geometric generation multiplies the boundary measure by $\rho>1$:
$$
|\Gamma_{\delta(t_{j+1})}| = \rho |\Gamma_{\delta(t_{j})}|.
$$
Writing $S_*:=|\Gamma_{\delta(t_0)}|$, the identity
$|\Gamma_{\delta(t_j)}|=S_*\rho^j$ shows that a stage which satisfies
$\varepsilon(t_{j_{\mathrm c}})|\Gamma_{\delta(t_{j_{\mathrm c}})}|\asymp1$,
obeys
\[
 j_{\mathrm c}
 =-\frac{\log\!\left(\varepsilon(t_{j_{\mathrm c}})S_*\right)}{\log\rho}
       +\mathcal O(1).
\]
This relation is implicit when the local coefficient evolves with time.
For a nonoverlapping equal-ratio self-similar construction satisfying the
open set condition, with $n$ pieces and similarity ratio $r<1$,
$\rho=nr^{d-1}=r^{-(D-d+1)}$ and $D=-\log n/\log r$; see
\cite{Hutchinson1981SelfSimilarity}. Hence,
\[
 j_{\mathrm c}
 =\frac{\log\!\left(\varepsilon(t_{j_{\mathrm c}})S_*\right)}
              {(D-d+1)\log r}+\mathcal O(1).
\]

As an example, we consider the planar Koch curve, where $d=2, n=4,$ and  $r= \frac{1}{3}$. As $\varepsilon(t_{j_{\mathrm c}})\to0$,
$$
 j_{\mathrm{c}}\sim - \frac{\log \varepsilon(t_{j_{\mathrm c}})}{\log (4/3)}.
$$

In terms of the approximation scale, suppose that the boundary approximations satisfy the sharper surface-growth law  $|\Gamma_\delta|\sim C_\Gamma\delta^{-\theta}$, where $\theta=D-(d-1)>0$; the existence of a fractal dimension alone does not imply this law.  If $\delta\sim\varepsilon^a$, then
\[
 \varepsilon|\Gamma_\delta|
          \sim C_\Gamma\varepsilon^{1-a\theta}.
\]
Thus, $a_{\mathrm c}=\theta^{-1}$ separates subcritical and supercritical couplings, with a finite nonzero total-mass limit at equality.
If only the two-sided comparability of the surface-growth law is known, the endpoint masses are bounded above and away from zero.  The operator interpretations require the additional hypotheses in Section~\ref{sec:operator-regimes}. 

\subsection{Numerical illustrations}\label{sec:numerical}

We use the unit-side Koch snowflake stages $\Omega_m$, $m=0,\ldots,6$, to illustrate Lemma~\ref{lem:uniform-prefractal-ground-state} and Theorem~\ref{thm:minkowski-recovery}. 
The initial domain is the equilateral triangle of side length one.  All quantities in the following are nondimensional.  The scale, perimeter, and limiting dimension are
\begin{equation}\label{eq:koch-geometry-numerical}
 \delta_m=3^{-m},\qquad |\Gamma_m|=3(4/3)^m,
 \qquad D=\frac{\log4}{\log3}=1.261859507\ldots.
\end{equation}
The area is also known exactly:
\begin{equation}\label{eq:koch-volume-numerical}
 |\Omega_m|=\frac{\sqrt3}{4}
       \left[\frac85-\frac35\left(\frac49\right)^m\right],
 \qquad |\Omega_m|\longrightarrow\frac{2\sqrt3}{5}.
\end{equation}
Indeed, generation $k\geq1$ adds $3\cdot4^{k-1}$ triangles of area
$(\sqrt3/4)9^{-k}$, and summing these additions gives
\eqref{eq:koch-volume-numerical}.  This variation of the area is relevant
to adjacent-stage dimension estimates, even though it does not affect
the logarithmic limit in Theorem~\ref{thm:minkowski-recovery}.

\subsubsection{Finite-element formulation}

For each fixed generation $m$, let $\mathcal T_{m,h}$ be a
conforming triangulation of the exact polygon $\Omega_m$, and define
the mesh size by
\[
 h:=\max_{T\in\mathcal T_{m,h}}\operatorname{diam}(T).
\]
We discretize $a_{\varepsilon,\delta_m}$ in the conforming
piecewise-linear space
\[
 V_{m,h}:=
 \left\{v\in C^0(\overline{\Omega_m}):
 v|_T\in\mathbb P_1(T)
 \text{ for every }T\in\mathcal T_{m,h}\right\},
\]
where \(\mathbb P_1(T)\) denotes the affine functions on \(T\).
Here, $h$ controls the numerical resolution on the fixed polygon,
whereas $\delta_m=3^{-m}$ specifies the geometric prefractal scale.
Let $\Lambda_{m,h}(\varepsilon)$ denote the exact discrete ground-state Robin eigenvalue. For a nodal basis $(\chi_p)$ of $V_{m,h}$,
the generalized eigenvalue problem is
\begin{equation}\label{eq:FEM-robin-matrices}
 (\mathbf K_m+\varepsilon\mathbf B_m)c
        =\Lambda_{m,h}(\varepsilon)\mathbf M_m c,
 \qquad c^{\mathsf T}\mathbf M_m c=1,
\end{equation}
where
\[
 (\mathbf K_m)_{pq}=\int_{\Omega_m}\nabla\chi_p\cdot\nabla\chi_q\,dx,
 \quad
 (\mathbf M_m)_{pq}=\int_{\Omega_m}\chi_p\chi_q\,dx,
 \quad
 (\mathbf B_m)_{pq}=\int_{\Gamma_m}\chi_p\chi_q\,d\sigma.
\]
The dependence of \(\mathbf K_m,\mathbf M_m,\mathbf B_m\) on \(h\) is suppressed in the notation. The boundary integral is assembled on the actual stage-$m$ polygon,
not on a smooth surrogate.  The exact assembly gives
\begin{equation}\label{eq:FEM-assembly-identities} 
      \mathbf K_m\mathbf1=0,\qquad \mathbf1^{\mathsf T}\mathbf M_m\mathbf1=|\Omega_m|,\qquad \mathbf1^{\mathsf T}\mathbf B_m\mathbf1=|\Gamma_m|. 
\end{equation}
Since constants belong to the conforming space, the discrete Rayleigh
principle gives the following useful bounds:
\begin{equation}\label{eq:FEM-consistency}
 \Lambda_m(\varepsilon)\leq\Lambda_{m,h}(\varepsilon)
                 \leq\frac{\varepsilon|\Gamma_m|}{|\Omega_m|}.
\end{equation}
The first inequality presupposes a conforming discretization of the
exact polygon and exact bilinear forms.

Moreover, set \(\tau:=\varepsilon|\Gamma_m|\). For \(0<\tau\leq\tau_0\),
combining the conforming finite-element bounds \eqref{eq:FEM-consistency} with Lemma \ref{lem:uniform-prefractal-ground-state} gives the following estimates:
\[
0\leq\Lambda_{m,h}(\varepsilon)-\Lambda_m(\varepsilon)
\leq\frac{\tau}{|\Omega_m|}-\Lambda_m(\varepsilon)
\leq C\tau^2.
\]
Consequently, the exact discrete eigenvalue also satisfies
\begin{equation} \label{computed-eigenvalue}
\Lambda_{m,h}(\varepsilon)
=\frac{\tau}{|\Omega_m|}\bigl(1+\mathcal O(\tau)\bigr),
\qquad \tau\to0,
\end{equation}
uniformly in the prefractal depth and in conforming spaces containing constants. For \(\varepsilon=\varepsilon_m\), one has \(\tau=\tau_m\).

\subsubsection{Dimension recovery and the finite-stage area correction}
\begin{figure}[htbp]
\centering
\includegraphics[width=0.98\textwidth]{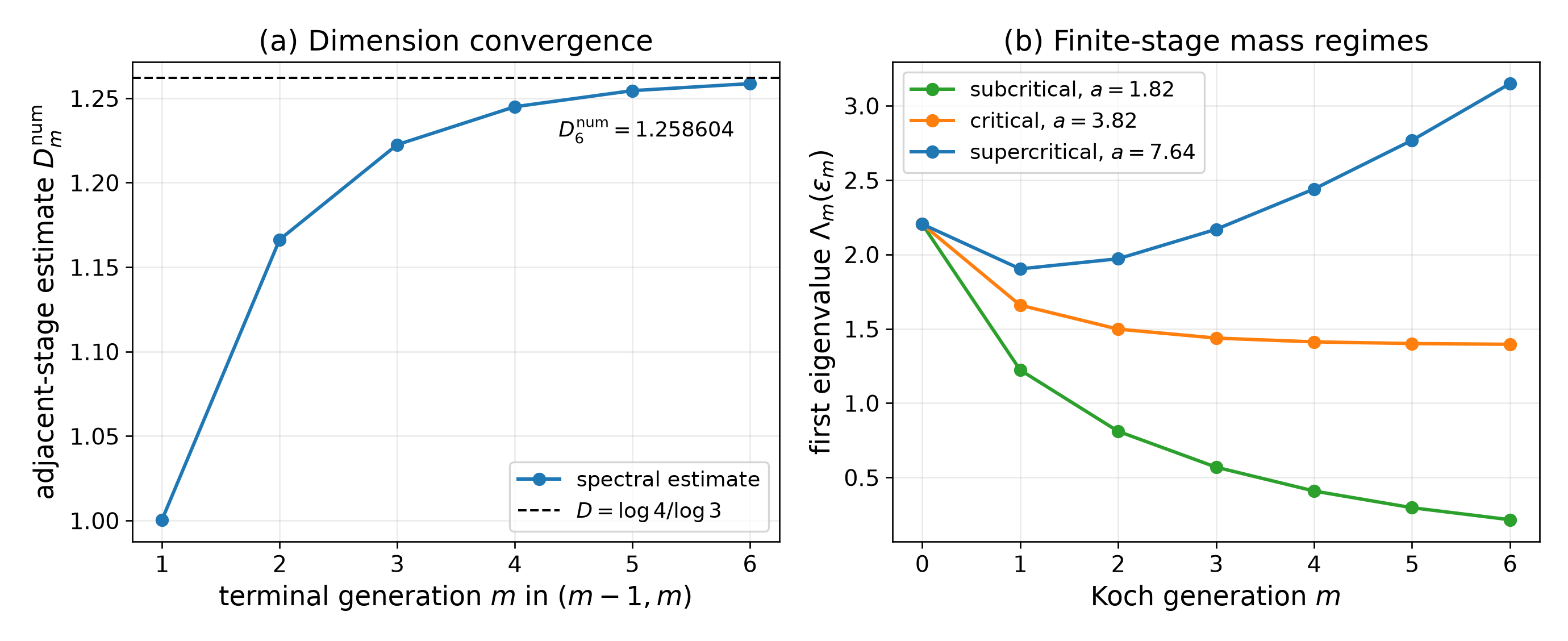}
\caption{Original finite-stage Koch computations, $m=0,\ldots,6$.
(a) The uncorrected estimator \eqref{eq:dimension-estimator-raw} along
$\varepsilon_m=10^{-2}\delta_m^{0.6}$; the dashed line is $D=\log4/\log3$.
Its terminal value is area-corrected in \eqref{eq:corrected-D6}; Table~\ref{tab:koch-refinement} gives a separate nested-mesh study of this experiment at $m=5,6$, including the corrected estimator. (b) Ground-state eigenvalues for the choices $a=20/11$, $a=a_{\mathrm c}$, and $a=2a_{\mathrm c}$, with $\varepsilon_m=\frac13\delta_m^{1/a}$.
The exact critical exponent is $a_{\mathrm c}=\log3/\log(4/3)$; the middle curve uses $a=a_{\mathrm c}$ and is exactly critical. The refinement table concerns panel~(a), not the curves in panel~(b). The parameter values in the legend are rounded to two decimal places.}
\label{fig:koch-numerics}
\end{figure}

The Koch family satisfies the exact perimeter ratio
\[
\frac{|\Gamma_m|}{|\Gamma_{m-1}|}=\frac43=3^{D-1}.
\]
This identity permits an adjacent-stage estimator whose error can be separated into area variation and the small-mass remainder. Along the subcritical diagonal
\begin{equation}\label{eq:koch-subcritical-diagonal}
 \varepsilon_m=10^{-2}\delta_m^{0.6},\qquad
 \tau_m=0.03\left(\frac4{3^{1.6}}\right)^m\longrightarrow0,
\end{equation}
we define the adjacent-stage estimator using continuum eigenvalues
\begin{equation}\label{eq:dimension-estimator-raw}
 D_m^{\mathrm{num}}
 :=1+\frac1{\log3}\log\left[
       \frac{\Lambda_m(\varepsilon_m)/\varepsilon_m}
            {\Lambda_{m-1}(\varepsilon_{m-1})/\varepsilon_{m-1}}
                        \right],\qquad m\geq 1.
\end{equation}
For computed data, inserting the corresponding discrete eigenvalues defines $D_{m,h}^{\mathrm{num}}$. Figure~\ref{fig:koch-numerics}(a) reports this \emph{uncorrected} estimator, including $D_{6,\mathrm{fig}}^{\mathrm{num}}=1.258604$.

For continuum data, set
$e_m:=|\Omega_m|\Lambda_m(\varepsilon_m)/\tau_m-1$.
Lemma~\ref{lem:uniform-prefractal-ground-state} gives
$e_m\leq0$ and $|e_m|\leq C\tau_m$ in the subcritical range.  From
\eqref{eq:koch-geometry-numerical}, it follows that
\begin{equation}\label{eq:raw-dimension-bias}
 D_m^{\mathrm{num}}
 =D-\frac{\log(|\Omega_m|/|\Omega_{m-1}|)}{\log3}
        +\frac{\log(1+e_m)-\log(1+e_{m-1})}{\log3}.
\end{equation}
The term $-\log(|\Omega_m|/|\Omega_{m-1}|)/\log3$ accounts for the change in domain area. 
In particular, replacing $\Lambda_m(\varepsilon_m)$ by its leading approximation $\tau_m/|\Omega_m|$ gives $D_1^{\mathrm{num}}=1$, since $|\Omega_1|/|\Omega_0|=4/3$.

The area-corrected estimator is therefore
\begin{align}\label{eq:dimension-estimator-volume}
 D_m^{\mathrm{vol}}
 &:=1+\frac1{\log3}\log\left[
       \frac{|\Omega_m|\Lambda_m(\varepsilon_m)/\varepsilon_m}
            {|\Omega_{m-1}|\Lambda_{m-1}(\varepsilon_{m-1})/
                                      \varepsilon_{m-1}}\right]\notag\\
 &=D_m^{\mathrm{num}}
       +\frac{\log(|\Omega_m|/|\Omega_{m-1}|)}{\log3}.
\end{align}
Since $\tau_m\to0$, both $|e_m|$ and $|e_{m-1}|$ are at most $1/2$ for all sufficiently large $m$. Applying $|\log(1+t)|\leq2|t|$ for $|t|\leq1/2$ to \eqref{eq:raw-dimension-bias} gives
\begin{equation}\label{eq:volume-estimator-error}
D_m^{\mathrm{vol}}
=D+\mathcal O(\tau_m+\tau_{m-1}),
\qquad m\to\infty.
\end{equation}
By \eqref{computed-eigenvalue}, the same small-mass estimate holds for $D_{m,h}^{\mathrm{vol}}$ evaluated with conforming discrete eigenvalues at the two generations. Thus, the correction removes the area contribution without another eigenvalue solve.

At generation six, the correction is
$\log(|\Omega_6|/|\Omega_5|)/\log3=0.003304050958\ldots$.
Applied to the rounded value in Figure~\ref{fig:koch-numerics}(a), it gives
\begin{equation}\label{eq:corrected-D6}
 D_{6,\mathrm{fig}}^{\mathrm{vol}}
:=D_{6,\mathrm{fig}}^{\mathrm{num}}+0.003304050958\ldots
 \simeq1.261908051.
\end{equation}
Its discrepancy from $D$ is about $4.85\times10^{-5}$, compared with
$3.26\times10^{-3}$ before correction.  This is an algebraic correction
of the original rounded data. The computations in
Table~\ref{tab:koch-refinement} below assess the mesh sensitivity for the
same two generations and the same coefficients.

\subsubsection{Mesh refinement}
\label{sec:koch-refinement}

We recomputed the ground-state eigenvalues on the fixed full polygons
$\Omega_5$ and $\Omega_6$, with \eqref{eq:koch-subcritical-diagonal},
on three nested conforming $P_1$ meshes.  The initial meshes use graded
triangular-lattice points and Delaunay triangulation; their boundary
edges coincide with the prescribed polygon edges.  Their minimum
angle is approximately $19.1^\circ$ and maximum element diameter $1/9$.
At each refinement, every triangle is divided into four by joining its
edge midpoints.  The polygon and Robin coefficient remain unchanged,
while $h=h_q=2^{-q}/9$ and the boundary edge length is $2^{-q}\delta_m$,
$q=0,1,2$.  The finest meshes contain $135\,552$ and $544\,896$
triangles, respectively.

We assemble exact local volume and boundary integrals and use a
generalized shift-invert Lanczos solver with shift $-1$ and tolerance
$10^{-12}$.  Eigenvalues are reevaluated by elementwise Rayleigh quotients. The assembly identities \eqref{eq:FEM-assembly-identities} are verified to floating-point precision.

Let $N_{m,q}$ denote the number of degrees of freedom and
$\Lambda_{m,q}:=\Lambda_{m,h_q}(\varepsilon_m)$.  The table reports
\[
 r_q:=\max_{m\in\{5,6\}}
       \frac{|\Lambda_{m,q}-\Lambda_{m,q-1}|}{\Lambda_{m,q}}
       \quad(q=1,2),
 \qquad
 D_{6,q}^{\mathrm{vol}}
 :=1+\frac1{\log3}\log\left[
       \frac{|\Omega_6|\Lambda_{6,q}/\varepsilon_6}
            {|\Omega_5|\Lambda_{5,q}/\varepsilon_5}\right].
\]
For nested conforming spaces, the exact discrete eigenvalues satisfy
\[
 \Lambda_m(\varepsilon_m)\leq\Lambda_{m,q+1}
       \leq\Lambda_{m,q}\leq\frac{\tau_m}{|\Omega_m|}.
\]
The computed values have this monotonicity.

\begin{table}[htbp]
\centering
\small
\setlength{\tabcolsep}{4pt}
\begin{tabular}{@{}r r r r r c r@{}}
\toprule
$q$ & $N_{5,q}$ & $N_{6,q}$ & $10^3\Lambda_{5,q}$ & $10^3\Lambda_{6,q}$
    & $r_q$ & $D_{6,q}^{\mathrm{vol}}$ \\
\midrule
0 & $5\,773$ & $23\,173$ & 6.801259135 & 4.674145955
  & -- & 1.261908544 \\
1 & $20\,017$ & $80\,401$ & 6.801244587 & 4.674140108
  & $2.139\times10^{-6}$ & 1.261909352 \\
2 & $73\,921$ & $297\,025$ & 6.801239809 & 4.674138256
  & $7.025\times10^{-7}$ & 1.261909631 \\
\bottomrule
\end{tabular}
\caption{Separate mesh-refinement computations for the final two
generations of the subcritical reconstruction experiment in
Figure~\ref{fig:koch-numerics}(a), with
$\varepsilon_m=10^{-2}3^{-0.6m}$.  Both fixed polygons are refined twice;
$N_{m,q}$ counts degrees of freedom and $r_q$ is the maximum relative
eigenvalue change.  The corrected estimate uses both eigenvalues at
level $q$. }
\label{tab:koch-refinement}
\end{table}

The maximum relative eigenvalue change on the final refinement is
$7.03\times10^{-7}$, and
\[
 |D_{6,2}^{\mathrm{vol}}-D_{6,1}^{\mathrm{vol}}|
       =2.79\times10^{-7},\qquad
 D_{6,2}^{\mathrm{vol}}=1.261909631\ldots,
 \qquad D_{6,2}^{\mathrm{vol}}-D\simeq5.01\times10^{-5}.
\]
For comparison, the constant-test upper bounds $\tau_m/|\Omega_m|$ are
$0.006802429150$ and $0.004674698200$ for $m=5,6$. The finest computed
eigenvalues lie below these upper bounds by relative amounts
$1.75\times10^{-4}$ and $1.20\times10^{-4}$, respectively, consistently
with the small-mass correction.

The finest \emph{uncorrected} value is
$D_{6,2}^{\mathrm{num}}=1.258605580\ldots$, close to, but distinct from,
the original value $1.258604$ in Figure~\ref{fig:koch-numerics}(a).
The last refinement change is much smaller than the remaining
discrepancy from $D$; the latter still contains the finite-Robin
perturbation contribution at fixed $m$ and $\varepsilon_m$.

\subsubsection{Interpretation and error analysis of the numerical results}
For Figure~\ref{fig:koch-numerics}(b), set $\varepsilon_m=\frac13\delta_m^{1/a}$. The total Robin boundary mass is given explicitly by
\begin{equation}
\label{eq:koch-three-mass-regimes}
\tau_m=\left(\frac43\,3^{-1/a}\right)^m,\qquad a_{\mathrm c}=\frac{\log3}{\log(4/3)}=3.818841679\ldots.
\end{equation}
The plotted choices $a=20/11$, $a=a_{\mathrm c}$, and $a=2a_{\mathrm c}$ are subcritical, critical, and supercritical, respectively. At exact criticality, $\tau_m=1$ at every stage. In the subcritical regime, $\Lambda_m(\varepsilon_m)\to0$. Under the corresponding hypotheses of Section~\ref{sec:operator-regimes}, the critical and supercritical problems have generalized Robin and Dirichlet-type limits, respectively.

We next assess how eigenvalue errors affect the dimension reconstruction. Let $\Lambda_m^{\mathrm{obs}}=\Lambda_m(\varepsilon_m)(1+\eta_m)$, where $\eta_m$ collects relative discretization, solver, and measurement errors. The induced error in \eqref{eq:dimension-estimator-volume} is
\begin{equation}\label{eq:dimension-error-propagation} D_m^{\mathrm{vol,obs}}-D_m^{\mathrm{vol}} =\frac{\log(1+\eta_m)-\log(1+\eta_{m-1})}{\log3}. \end{equation}
For $|\eta_m|,|\eta_{m-1}|\leq1/2$, its absolute value is at most $2(|\eta_m|+|\eta_{m-1}|)/\log3$. This error is separate from the finite-mass perturbation term in \eqref{eq:volume-estimator-error}. A common multiplicative error at adjacent stages cancels in the ratio, so Table~\ref{tab:koch-refinement} reports both the individual eigenvalues and the dimension estimate. The refinement differences provide evidence of discretization stability.

\section{Concluding remarks}\label{sec:conclusion}

In this work, the total Robin boundary mass $\tau_m=\varepsilon_m|\Gamma_m|$ links the local Robin coefficient to the spectral effect of increasingly rough prefractal boundaries.
A common normalized Poincar\'e-trace estimate gives eigenvalue bounds and, as $\tau_m\to0$, a first-order ground-state expansion with constants independent of the approximation level $m$.
It also identifies the rank-one resolvent singularity when $\tau_m\to0$.  Under the boundary-resolution hypothesis and the condition on the approximation scales, ground-state measurements with known Robin coefficients recover the lower and upper Minkowski dimensions for every sequence with bounded total mass. 

The operator interpretation requires more information than the total boundary mass.  Under the measure-convergence and compact Mosco hypotheses specified in Section~\ref{sec:operator-regimes}, $\tau_m\to0$ gives a Neumann limit, whereas $\tau_m\to\tau_*\in(0,\infty)$ gives a generalized Robin limit. For $\tau_m\to\infty$, the Dirichlet-type limit additionally requires the stated zero-trace and recovery conditions.
When the domain volumes are known, subcritical ground-state measurements determine the total Robin boundary mass to leading order through $\tau_m\sim|\Omega_m|\Lambda_m(\varepsilon_m)$. These leading-order data do not separate $|\Gamma_m|$ from an unknown $\varepsilon_m$. Independent calibration or additional geometric or modal information is necessary for joint reconstruction.

The layer-potential argument in Appendix~\ref{sec:proofs} gives the fixed-domain eigenpair coefficients and a remainder controlled by the spectral gap.  
For the Koch approximations, the exact area formula allows us to correct the adjacent-stage dimension estimator for the change in domain area between successive stages.
Extending the high-frequency estimates uniformly over specified prefractal families would connect the fixed-domain theory to deeper boundary approximations. A further question is whether measurements from several eigenmodes can separate boundary geometry from an unknown Robin coefficient.

\section*{Statements and declarations}

\paragraph{Funding.}
Yat Tin Chow was supported by the National Science Foundation under Grant No.~DMS-2409903 and by the Office of Naval Research under Award No.~N00014-24-1-2661.

\paragraph{Data availability.}
No external datasets were used in this study. 

\paragraph{Code availability.}
The code used to produce the numerical experiments is available from the corresponding author upon reasonable request.

\paragraph{Use of generative AI.}
During the preparation of this manuscript, the authors used ChatGPT (OpenAI) for language editing and assistance with the development of code for the numerical experiments. All AI-assisted text, code, and numerical outputs were independently reviewed and verified by the authors. The authors take full responsibility for the content of the article.

\appendix
\section{Proofs of the asymptotic formulas} \label{sec:proofs}

We retain the notation $\lambda_0^j$, $u_0^j$, $\lambda_\varepsilon^j$, and
$u_\varepsilon^j$ introduced in Section~\ref{sec:fixed-domain}.  The eigenfunctions of the Neumann problem are chosen to be real and orthonormal in $L^2(\Omega)$. The fixed-domain expansions are proved below by the layer-potential and Gohberg--Sigal method of \cite[Chapters~1--3]{surv153}; see also \cite[Sections~2.8--2.9]{surv235}.  The eigenvalue coefficients are contour residues, and the eigenfunction correction is obtained from a finite-rank factorization.
The remaining subsections treat the zero mode, multiple eigenvalues, and the spectral-gap dependence of the remainder.

The free-space layer-potential characterization is stated in the following for a bounded smooth connected domain with a connected exterior.  No constants uniform in a prefractal approximation are asserted in this appendix.  The spectral
transfer operators used below also give the same perturbation formulas on
smooth connected domains with holes; the qualification concerning the
free-space characteristic operator is explained at the end of the first
subsection.

The spectral transfer-operator argument below also proves the fixed-mode formulas under the $C^{1,\alpha}$ boundary assumption of Section~\ref{sec:fixed-domain}; the sharp boundary restriction estimate \eqref{tataru} is used only under the $C^\infty$ assumption.

\subsection{Robin eigenvalue characterization}

Let $G^\omega$ denote the outgoing fundamental solution with the convention
$(\Delta+\omega^2)G^\omega=\delta_0$.  Thus, for positive real $\omega$,
\[
G^\omega(x)=
\begin{cases}
-\dfrac{i}{4}H_0^{(1)}(\omega|x|),&d=2,\\[1mm]
-\dfrac{e^{i\omega|x|}}{4\pi|x|},&d=3.
\end{cases}
\]
Define the single- and double-layer potentials by
\[
S_\Omega^\omega[\varphi](x)
 =\int_{\partial\Omega}G^\omega(x-y)\varphi(y)\,d\sigma(y),
\qquad
D_\Omega^\omega[\varphi](x)
 =\int_{\partial\Omega}\frac{\partial G^\omega(x-y)}{\partial\nu(y)}
           \varphi(y)\,d\sigma(y).
\]
As usual, $S_\Omega^\omega$ also denotes the continuous boundary trace of the
single layer when its range is a boundary space.  The interior and exterior
traces are denoted by $-$ and $+$, respectively.  With the outward normal
convention of \cite[Section~2.3]{surv153},
\begin{equation}\label{eq:jumpD}
D_\Omega^\omega[\varphi]\big|_\pm
 =\left(\mp\frac12I+K_\Omega^\omega\right)[\varphi],
\qquad
\frac{\partial S_\Omega^\omega[\varphi]}{\partial\nu}\Big|_\pm
 =\left(\pm\frac12I+K_\Omega^{\omega,*}\right)[\varphi].
\end{equation}
Here, $K_\Omega^\omega$ differentiates the kernel in the source normal, whereas
$K_\Omega^{\omega,*}$ differentiates it in the target normal.  For a complex
outgoing kernel, this last notation should not be confused with the
Hilbert-space adjoint at the same complex frequency.

The boundary unknown in this appendix is the trace $\phi=u|_{\partial\Omega}$,
not a single-layer density.  Green's layer representation is
\begin{equation}\label{eq:GreenRep}
u=D_\Omega^\omega[u|_{\partial\Omega}]
       -S_\Omega^\omega[\partial_\nu u]
       \quad\hbox{in }\Omega.
\end{equation}
For the Robin condition, this becomes
\begin{equation}\label{eq:RobinRep}
u=D_\Omega^\omega[\phi]+\varepsilon S_\Omega^\omega[\phi].
\end{equation}
Taking its interior trace gives
\begin{equation}\label{eq:Adef}
 A(\omega,\varepsilon)[\phi]=0,
 \qquad
 A(\omega,\varepsilon):=-\frac12I+K_\Omega^\omega
                             +\varepsilon S_\Omega^\omega.
\end{equation}
In particular, the operator here involves $K_\Omega^\omega$, not
$K_\Omega^{\omega,*}$.  A formulation using the latter operator is possible,
but its unknown and reconstruction formulas are different.

\begin{lemma}[Robin eigenvalue characterization]\label{lem:char}
Assume that $\mathbb R^d\setminus\overline\Omega$ is connected.  For
$\omega>0$ and real $\varepsilon$, $\omega^2$ is a Robin eigenvalue if and only
if $\omega$ is a positive real characteristic value of
$A(\cdot,\varepsilon):L^2(\partial\Omega)\to L^2(\partial\Omega)$.
The correspondence is given by \eqref{eq:RobinRep}, with
$\phi=u|_{\partial\Omega}$.
\end{lemma}
\begin{proof}
A nonzero Robin eigenfunction has nonzero boundary trace, since otherwise
both Cauchy traces would vanish.  Equations~\eqref{eq:GreenRep}--\eqref{eq:Adef}
prove one implication.  Conversely, let $A(\omega,\varepsilon)\phi=0$ with
$\phi\ne0$, and set
$U=D_\Omega^\omega\phi+\varepsilon S_\Omega^\omega\phi$ on both sides of the
boundary.  Its exterior trace is $A(\omega,\varepsilon)\phi=0$.
Exterior Dirichlet uniqueness and the radiation condition therefore give
$U=0$ in the exterior.  The double layer has trace jump
$U|_--U|_+=\phi$, and its normal derivative has no jump; the single-layer
normal derivative has jump $\varepsilon\phi$.  Hence,
\[
U|_-=\phi,\qquad \partial_\nu U|_-=-\varepsilon\phi.
\]
The interior field satisfies the Robin condition and is not zero, because
its trace is $\phi\ne0$.
\end{proof}

For smooth $\partial\Omega$, $K_\Omega^\omega$ and
$S_\Omega^\omega:L^2(\partial\Omega)\to L^2(\partial\Omega)$ are compact.
Consequently $A$ is Fredholm of index zero.  It is holomorphic near every
positive frequency in both dimensions.  In dimension two, the logarithmic
branch at $\omega=0$ is excluded from this local assertion.
For the normalized branch, we write, as before,
\begin{equation}\label{eq:root}
 A(\omega_\varepsilon^l,\varepsilon)\phi_\varepsilon^l=0,
 \qquad \omega_\varepsilon^l=\sqrt{\lambda_\varepsilon^l},
\end{equation}
and
\begin{equation}\label{eq:reconstruct}
 u_\varepsilon^l
 =D_\Omega^{\omega_\varepsilon^l}[\phi_\varepsilon^l]
       +\varepsilon S_\Omega^{\omega_\varepsilon^l}[\phi_\varepsilon^l].
\end{equation}

If the complement has bounded components, their Dirichlet eigenfrequencies
may give additional characteristic values of this particular free-space
operator. We must not identify these with Robin eigenvalues in $\Omega$.
The boundary Neumann transfer operator in \eqref{eq:N-spectral} is intrinsic
to $\Omega$ and has only the Neumann poles of $\Omega$.  The regularized
block pencil in \eqref{eq:regularized-pencil} therefore provides an
alternative characterization that avoids these exterior artifacts.

\subsection{Gohberg--Sigal contour formula and normalization of the pole}

Fix $l\geq1$ with $\lambda_0^l$ being simple, and put
\[
\omega_0^l=\sqrt{\lambda_0^l},\qquad
A_0(\omega):=A(\omega,0)=-\frac12I+K_\Omega^\omega,
\qquad \phi_0^l=u_0^l|_{\partial\Omega}.
\]
The Neumann eigenvalue characterization and pole theorem imply that
$\omega_0^l$ is a simple characteristic value; see
\cite[Propositions~2.5--2.6 and Lemma~3.6]{surv153}.
Choose a small disk $V$ centered at $\omega_0^l$ whose closure meets no other
characteristic value and, in two dimensions, avoids the logarithmic cut.
The generalized Rouch\'e theorem applies whenever
\[
 |\varepsilon|\sup_{\omega\in\partial V}
       \|A_0(\omega)^{-1}S_\Omega^\omega\|<1.
\]
It gives exactly one characteristic value $\omega_\varepsilon^l$ in $V$,
counted with multiplicity.  The generalized argument principle then gives
\begin{equation}\label{eq:GSA}
 \omega_\varepsilon^l-\omega_0^l
 =\frac{1}{2\pi i}\operatorname{tr}
    \int_{\partial V}(\omega-\omega_0^l)
       A(\omega,\varepsilon)^{-1}
       \partial_\omega A(\omega,\varepsilon)\,d\omega.
\end{equation}
The contour expression also proves holomorphic dependence on
$\varepsilon$.  For real $\varepsilon$, the nearby eigenvalue is real.

The residue of $A_0^{-1}$ is not an unnormalized orthogonal projection.
With the inner product linear in its second argument, and with
$\phi_l^*\ne0$ in $\ker A_0(\omega_0^l)^*$, its precise form is
\begin{equation}\label{pole-pencil}
 A_0(\omega)^{-1}q
 =\frac{(\phi_l^*,q)_{\partial\Omega}\,\phi_0^l}
 {(\omega-\omega_0^l)
  (\phi_l^*,\partial_\omega A_0(\omega_0^l)\phi_0^l)_{\partial\Omega}}
       +R(\omega)q.
\end{equation}
Here, $R$ is holomorphic and the denominator is nonzero by simplicity of
the characteristic value.  The adjoint in this formula is the actual
Hilbert-space adjoint of $A_0(\omega_0^l)$.

On $\partial V$, the following expansion converges in operator norm:
\begin{align}\label{eq:NeumannSeriesA}
 A(\omega,\varepsilon)^{-1}
 &=\bigl(I+\varepsilon A_0(\omega)^{-1}S_\Omega^\omega\bigr)^{-1}
                       A_0(\omega)^{-1}\notag\\
 &=\sum_{p=0}^{\infty}(-\varepsilon)^p
       \bigl(A_0(\omega)^{-1}S_\Omega^\omega\bigr)^p A_0(\omega)^{-1}.
\end{align}
Substitution in \eqref{eq:GSA}, cyclicity of contour-residue traces, and
integration by parts give the Ammari--Kang--Lee (AKL) expansion
\begin{equation}\label{eq:AKL-general-moment}
 F(\omega_\varepsilon^l)-F(\omega_0^l)
 =\sum_{p=1}^{\infty}\frac{(-\varepsilon)^p}{2\pi i\,p}
       \operatorname{tr}\int_{\partial V}F'(\omega)
          \bigl(A_0(\omega)^{-1}S_\Omega^\omega\bigr)^p\,d\omega
\end{equation}
for every scalar holomorphic $F$ on a neighborhood of $\overline V$ and
$\varepsilon$ sufficiently small.  In particular, the first-order
\emph{frequency} term is
\[
 -\frac{\varepsilon}{2\pi i}\operatorname{tr}
          \int_{\partial V}A_0(\omega)^{-1}S_\Omega^\omega\,d\omega.
\]
For the \emph{eigenvalue} $\lambda_\varepsilon^l=(\omega_\varepsilon^l)^2$,
we take $F(\omega)=\omega^2$ and hence, $F'(\omega)=2\omega$.
This distinction fixes both the sign and the frequency factor.

For completeness, let $M=A_0^{-1}S_\Omega^\omega$.  The identity
$A=A_0(I+\varepsilon M)$ implies
\[
 A^{-1}\partial_\omega A
 =(I+\varepsilon M)^{-1}A_0^{-1}(\partial_\omega A_0)(I+\varepsilon M)
      +(I+\varepsilon M)^{-1}\varepsilon\partial_\omega M.
\]
In the contour trace the first term has the same contribution as
$A_0^{-1}\partial_\omega A_0$.  Expanding the second term gives
$\sum_{p\geq1}(-1)^{p-1}\varepsilon^pM^{p-1}\partial_\omega M$.
Cyclicity replaces its contour trace by
$p^{-1}\partial_\omega(M^p)$, and integration by parts yields
\eqref{eq:AKL-general-moment}, exactly as in
\cite[Theorem~3.9]{surv153}.

Throughout, the trace is taken \emph{after the contour integration}.
The integral is a finite-rank residue; the operator-valued integrand need not be trace class. The spectral representation below allows us to evaluate these contour residues.  

\subsection{The Neumann-to-Dirichlet operator and its spectral residues}

Write $X=L^2(\partial\Omega)$ and
\[
 f_j:=u_0^j|_{\partial\Omega},\qquad
 a_{lj}=(f_l,f_j)_{\partial\Omega}
       =\int_{\partial\Omega}u_0^lu_0^j\,d\sigma,\qquad a_l:=a_{ll}>0.
\]
For $f,g\in X$, the notation $f\otimes g$ denotes
$q\mapsto(g,q)_{\partial\Omega}f$.  The rank-one operator
$\Pi_l:=f_l\otimes f_l$ is generally not a projection:
$\Pi_l^2=a_l\Pi_l$ and $\operatorname{tr}\Pi_l=a_l$.

On a neighborhood of $\lambda_0^l>0$, choose the branch
$\omega=\sqrt z$ with $\sqrt{\lambda_0^l}=\omega_0^l$.  Define
\begin{equation}\label{eq:boundary-transfer}
 \mathcal N(z):=A_0(\sqrt z)^{-1}S_\Omega^{\sqrt z},
 \qquad
 \mathcal P(z):=D_\Omega^{\sqrt z}\mathcal N(z)-S_\Omega^{\sqrt z}.
\end{equation}
The second expression is an interior potential, mapping $X$ to
$H^1(\Omega)$. The unperturbed layer identities in
\cite[Propositions~2.7--2.8 and formula~(2.30)]{surv153} identify
$\mathcal N$ as the boundary Neumann-to-Dirichlet operator and
$\mathcal P$ as its interior Neumann potential. Thus,
\begin{equation}\label{eq:N-spectral}
 \mathcal N(z)q
 =\sum_{j=0}^{\infty}
       \frac{(f_j,q)_{\partial\Omega}f_j}{\lambda_0^j-z},
 \qquad
 \mathcal P(z)q
 =\sum_{j=0}^{\infty}
       \frac{(f_j,q)_{\partial\Omega}u_0^j}{\lambda_0^j-z}.
\end{equation}
These are identities for the \emph{unperturbed} transfer operators.
They imply, in particular, that the assembled operators are
single-valued meromorphic functions of $z$, including near $z=0$.
This assertion is about $\mathcal N$ and $\mathcal P$, not about each
individual outgoing layer at $\sqrt z=0$.

Here is a norm justification of \eqref{eq:N-spectral} and of the series
used in the following.  Let $\mathcal L_0$ be the unperturbed Neumann Laplacian and
$\gamma:H^1(\Omega)\to X$ its trace map.  The operator
\[
 \mathcal T:=\gamma(1+\mathcal L_0)^{-1/2}:L^2(\Omega)\to X
\]
is compact.  This follows from the identification
$D((1+\mathcal L_0)^{1/2})=H^1(\Omega)$ and the compact embedding of the
trace space $H^{1/2}(\partial\Omega)$ in $X$.
Since $\mathcal T u_0^j=f_j/\sqrt{1+\lambda_0^j}$, Parseval's identity gives
\begin{equation}\label{eq:weighted-summability-A}
 \sum_{j=0}^{\infty}
       \frac{|(f_j,q)_{\partial\Omega}|^2}{1+\lambda_0^j}
 =\|\mathcal T^*q\|_{L^2(\Omega)}^2
 \leq\|\mathcal T\|^2\|q\|_X^2.
\end{equation}
Furthermore,
\begin{align}\label{eq:N-factorization-A}
 \mathcal N(z)&=\mathcal T(1+\mathcal L_0)
                    (\mathcal L_0-z)^{-1}\mathcal T^*,\notag\\
 \mathcal P(z)&=(1+\mathcal L_0)^{-1/2}
                    (1+\mathcal L_0)(\mathcal L_0-z)^{-1}\mathcal T^*.
\end{align}
For $z$ in a compact subset of the complement of the Neumann spectrum,
the middle diagonal multiplier is uniformly bounded.  Finite spectral
truncations in \eqref{eq:N-factorization-A} therefore converge in operator
norm, using the compactness of $\mathcal T$ or $\mathcal T^*$.
For the second identity, the target norm is $H^1(\Omega)$, with its
Neumann spectral energy norm. Thus, the two series in
\eqref{eq:N-spectral} have precisely the required operator interpretations.

Near the fixed simple eigenvalue, separate the pole:
\begin{equation}\label{eq:N-pole-A}
 \mathcal N(z)=-\frac{\Pi_l}{z-\lambda_0^l}+\mathcal N_l(z),
 \qquad
 \mathcal P(z)=-\frac{u_0^l\otimes f_l}{z-\lambda_0^l}
                                     +\mathcal P_l(z).
\end{equation}
Here, both regular parts are holomorphic on a disk containing no other
Neumann eigenvalue, and
\begin{equation}\label{eq:N-regular-A}
 \mathcal N_l(\lambda_0^l)q
 =\sum_{j\ne l}\frac{(f_j,q)_{\partial\Omega}f_j}{\lambda_0^j-\lambda_0^l},
 \qquad
 \mathcal P_l(\lambda_0^l)q
 =\sum_{j\ne l}\frac{(f_j,q)_{\partial\Omega}u_0^j}
                              {\lambda_0^j-\lambda_0^l}.
\end{equation}
In particular,
\begin{equation}\label{eq:absolute-B-A}
 \sum_{j\ne l}\frac{|a_{lj}|^2}{|\lambda_0^j-\lambda_0^l|}<\infty,
 \qquad
 \sum_{j\ne l}\frac{(1+\lambda_0^j)|a_{lj}|^2}
                         {|\lambda_0^j-\lambda_0^l|^2}<\infty.
\end{equation}
Indeed, for fixed $l$, each weight is bounded by a constant depending on
$l$ times $(1+\lambda_0^j)^{-1}$; apply
\eqref{eq:weighted-summability-A} with $q=f_l$.
This proves both absolute convergence of $B_l$ and $H^1$ convergence of
the eigenfunction correction. 

\subsection{First two terms in the asymptotic expansions of the eigenvalues}

For $z$ near $\lambda_0^l>0$, the original operator factorizes as
\begin{equation}\label{eq:factor-A-N}
 A(\sqrt z,\varepsilon)
 =A_0(\sqrt z)\bigl(I+\varepsilon\mathcal N(z)\bigr).
\end{equation}
With $F(\omega)=\omega^2$, the change of variables $z=\omega^2$ in
\eqref{eq:AKL-general-moment} gives
\begin{equation}\label{eq:lambda-contour-series}
 \lambda_\varepsilon^l-\lambda_0^l
 =\sum_{p=1}^{\infty}\frac{(-\varepsilon)^p}{2\pi i\,p}
       \operatorname{tr}\int_{\partial W_l}\mathcal N(z)^p\,dz,
\end{equation}
where $W_l$ is a sufficiently small disk centered at $\lambda_0^l$.
Only the principal parts at $z=\lambda_0^l$ contribute.
For $p=1$, \eqref{eq:N-pole-A} gives
\[
 \frac1{2\pi i}\operatorname{tr}
        \int_{\partial W_l}\mathcal N(z)\,dz=-a_l.
\]
Consequently, the first coefficient is $a_l=\int_{\partial\Omega}(u_0^l)^2$.
For $p=2$, write $\zeta=z-\lambda_0^l$ and expand
\[
 \mathcal N(z)^2
 =\frac{\Pi_l^2}{\zeta^2}
  -\frac{\Pi_l\mathcal N_l(\lambda_0^l)
             +\mathcal N_l(\lambda_0^l)\Pi_l}{\zeta}
  +\hbox{a holomorphic operator}.
\]
The double pole has zero contour integral. Therefore,
\[
 \frac1{2\pi i}\operatorname{tr}\int_{\partial W_l}\mathcal N(z)^2\,dz
 =-2(f_l,\mathcal N_l(\lambda_0^l)f_l)_{\partial\Omega}.
\]
Substitution into \eqref{eq:lambda-contour-series} yields
\begin{equation}\label{eq:lambda-residue-evaluation}
 \lambda_\varepsilon^l
 =\lambda_0^l+\varepsilon a_l
       -\varepsilon^2(f_l,\mathcal N_l(\lambda_0^l)f_l)_{\partial\Omega}
       +\mathcal O_l(\varepsilon^3).
\end{equation}
By \eqref{eq:N-regular-A} and \eqref{eq:absolute-B-A},
\begin{equation}\label{eq:B-residue-A}
 -(f_l,\mathcal N_l(\lambda_0^l)f_l)_{\partial\Omega}
 =\sum_{j\ne l}\frac{|a_{lj}|^2}{\lambda_0^l-\lambda_0^j}=B_l.
\end{equation}
This proves \eqref{asymp:eigenvaluel} directly from the first two AKL
contour residues.  Expanding the square root also gives
\begin{equation}\label{eq:frequency-two-terms-A}
 \omega_\varepsilon^l
 =\omega_0^l+\frac{a_l}{2\omega_0^l}\varepsilon
 +\left(\frac{B_l}{2\omega_0^l}
           -\frac{a_l^2}{8(\omega_0^l)^3}\right)\varepsilon^2
 +\mathcal O_l(\varepsilon^3).
\end{equation}
The eigenvalue and frequency second coefficients must therefore not be
identified with each other.

\subsection{Finite-rank factorization and the eigenfunction expansion}

To reconstruct the eigenfunction and control the remainder, set
\[
 Q_l(z,\varepsilon):=I+\varepsilon\mathcal N_l(z),
 \qquad
 \psi_l(z,\varepsilon):=Q_l(z,\varepsilon)^{-1}f_l.
\]
For $\varepsilon$ sufficiently small, $Q_l$ is holomorphic and invertible on
$\overline W_l$.  Formula~\eqref{eq:N-pole-A} gives the exact factorization
\begin{equation}\label{eq:rank-one-factor-A}
 I+\varepsilon\mathcal N(z)
 =Q_l(z,\varepsilon)
       \left(I-\frac{\varepsilon}{z-\lambda_0^l}
                       \psi_l(z,\varepsilon)\otimes f_l\right).
\end{equation}
The last factor is a rank-one perturbation of the identity.  Its
characteristic equation is
\begin{equation}\label{eq:scalar-characteristic-A}
 z-\lambda_0^l
   =\varepsilon(f_l,\psi_l(z,\varepsilon))_{\partial\Omega}.
\end{equation}

An analytic pencil that includes the removable unperturbed point is
\begin{equation}\label{eq:regularized-pencil}
 \mathfrak A_l(z,\varepsilon)
 =\begin{pmatrix}
 Q_l(z,\varepsilon)&-f_l\\
 -\varepsilon(f_l,\cdot)_{\partial\Omega}&z-\lambda_0^l
 \end{pmatrix}:X\oplus\mathbb C\longrightarrow X\oplus\mathbb C.
\end{equation}
Block Gaussian elimination reduces it, by holomorphic invertible factors,
to $\operatorname{diag}(Q_l,s_l)$, where
\[
 s_l(z,\varepsilon)=z-\lambda_0^l
              -\varepsilon(f_l,Q_l(z,\varepsilon)^{-1}f_l)_{\partial\Omega}.
\]
At $\varepsilon=0$ this scalar factor has a simple zero.  Rouch\'e's theorem
and its contour moment give a unique holomorphic zero
$z=\lambda_\varepsilon^l$ for small $\varepsilon$.
For real $\varepsilon$, the real spectral data in \eqref{eq:N-spectral}
give $s_l(\overline z,\varepsilon)=\overline{s_l(z,\varepsilon)}$;
uniqueness of the zero therefore makes $\lambda_\varepsilon^l$ real.
Since $a_l>0$, \eqref{eq:scalar-characteristic-A} also implies
$\lambda_\varepsilon^l-\lambda_0^l=\varepsilon a_l+\mathcal O_l(\varepsilon^2)$;
thus the shifted root does not coincide with the old pole when
$\varepsilon\ne0$ is sufficiently small.

At that root, write
$\psi_\varepsilon^l:=\psi_l(\lambda_\varepsilon^l,\varepsilon)$ and define
\begin{equation}\label{eq:interior-unnormalized-A}
 \widetilde u_\varepsilon^l
 :=-\varepsilon\mathcal P(\lambda_\varepsilon^l)\psi_\varepsilon^l.
\end{equation}
Since $(I+\varepsilon\mathcal N(\lambda_\varepsilon^l))
\psi_\varepsilon^l=0$, the unperturbed layer identity
\eqref{eq:boundary-transfer} gives
\[
 \widetilde u_\varepsilon^l
 =D_\Omega^{\omega_\varepsilon^l}[\psi_\varepsilon^l]
        +\varepsilon S_\Omega^{\omega_\varepsilon^l}[\psi_\varepsilon^l],
 \qquad \gamma\widetilde u_\varepsilon^l=\psi_\varepsilon^l.
\]
Thus, this is the Robin eigenfunction reconstructed by the original layer potentials.
Its component in the $u_0^l$ direction is exactly one: by
\eqref{eq:N-pole-A} and \eqref{eq:scalar-characteristic-A},
\begin{equation}\label{eq:exact-eigenfunction-A}
 \widetilde u_\varepsilon^l
 =u_0^l-\varepsilon\mathcal P_l(\lambda_\varepsilon^l)
                                        \psi_\varepsilon^l,
 \qquad
 (u_0^l,\widetilde u_\varepsilon^l)_{L^2(\Omega)}=1.
\end{equation}
The Neumann series for $Q_l^{-1}$ and holomorphy of the regular parts give
\[
 \psi_\varepsilon^l
   =f_l-\varepsilon\mathcal N_l(\lambda_0^l)f_l
                          +\mathcal O_X(\varepsilon^2),
\]
and
\begin{equation}\label{eq:eigenfunction-direct-A}
 \widetilde u_\varepsilon^l
 =u_0^l-\varepsilon\mathcal P_l(\lambda_0^l)f_l
                         +\mathcal O_{H^1(\Omega)}(\varepsilon^2)
 =u_0^l+\varepsilon\sum_{j\ne l}
      \frac{a_{lj}}{\lambda_0^l-\lambda_0^j}u_0^j
                         +\mathcal O_{H^1(\Omega)}(\varepsilon^2).
\end{equation}
The series converges in $H^1(\Omega)$ by \eqref{eq:absolute-B-A}.
Because the correction in \eqref{eq:exact-eigenfunction-A} is orthogonal to
$u_0^l$,
\[
 \|\widetilde u_\varepsilon^l\|_{L^2(\Omega)}^2
 =1+\varepsilon^2
   \|\mathcal P_l(\lambda_\varepsilon^l)\psi_\varepsilon^l\|_{L^2(\Omega)}^2
 \quad\hbox{for real }\varepsilon.
\]
Dividing by this norm changes the expansion only at order $\varepsilon^2$.
Consequently, with
$(u_0^l,u_\varepsilon^l)_{L^2(\Omega)}>0$, \eqref{asymp:eigenfunction}
holds in $H^1(\Omega)$, and its normalized boundary root is
$\phi_\varepsilon^l=\psi_\varepsilon^l/
\|\widetilde u_\varepsilon^l\|_{L^2(\Omega)}$.

The scalar factor also gives a direct remainder check on
\eqref{eq:lambda-residue-evaluation}.  Uniformly for $z\in W_l$,
\[
 (f_l,Q_l(z,\varepsilon)^{-1}f_l)
 =a_l-\varepsilon(f_l,\mathcal N_l(\lambda_0^l)f_l)
          +\mathcal O_l(\varepsilon|z-\lambda_0^l|+\varepsilon^2).
\]
At the root, $|z-\lambda_0^l|=\mathcal O_l(\varepsilon)$, so multiplication by $\varepsilon$ in \eqref{eq:scalar-characteristic-A} gives an $\mathcal O_l(\varepsilon^3)$ remainder. 

\subsection{The zero Neumann eigenvalue}

The mode $\lambda_0^0=0$ is simple, with
$u_0^0=|\Omega|^{-1/2}$ and $f_0=|\Omega|^{-1/2}$ on $\partial\Omega$.
Its frequency is of order $\sqrt\varepsilon$, so it cannot be treated as
an analytic function of $\varepsilon$ in the variable $\omega$ at zero.
In dimension two the separate outgoing layers also have logarithmic
singularities there.

We instead use the meromorphic transfer operators
\eqref{eq:N-spectral}, whose regular parts after subtraction of the zero
mode are holomorphic at $z=0$.  The block pencil
\eqref{eq:regularized-pencil}, with $l=0$, is holomorphic in $(z,\varepsilon)$
and has a simple characteristic value at $(0,0)$.  Its scalar Schur
complement is
\[
 z-\varepsilon(f_0,(I+\varepsilon\mathcal N_0(z))^{-1}f_0).
\]
The same finite-rank construction and interior reconstruction therefore
apply without a zero-frequency branch ambiguity.  Since
\[
 a_0=\frac{|\partial\Omega|}{|\Omega|},\qquad
 a_{0j}=\frac1{\sqrt{|\Omega|}}\int_{\partial\Omega}u_0^j\,d\sigma,
\]
they give \eqref{asym_lambda-fixed} and \eqref{asymp:eigenfunction0}, namely
\begin{align}\label{eq:zero-mode-appendix}
 \lambda_\varepsilon^0
 &=\frac{|\partial\Omega|}{|\Omega|}\varepsilon
 -\frac{\varepsilon^2}{|\Omega|}\sum_{j=1}^{\infty}
       \frac{\left(\int_{\partial\Omega}u_0^j\,d\sigma\right)^2}{\lambda_0^j}
       +\mathcal O(\varepsilon^3),\notag\\
 u_\varepsilon^0
 &=\frac1{\sqrt{|\Omega|}}
 -\frac{\varepsilon}{\sqrt{|\Omega|}}
       \sum_{j=1}^{\infty}\frac{\int_{\partial\Omega}u_0^j\,d\sigma}{\lambda_0^j}
                                     u_0^j
       +\mathcal O_{H^1(\Omega)}(\varepsilon^2).
\end{align}
Both convergence assertions follow from \eqref{eq:weighted-summability-A}.
The block formulation is a pole subtraction for assembled layer-potential
operators; it does not assume that $A(\sqrt z,\varepsilon)$ itself is
holomorphic at $z=0$.

\subsection{Splitting of a multiple eigenvalue}

Let $\lambda_0^l$ have multiplicity $m$, and retain the orthonormal basis
$u_0^{l,1},\ldots,u_0^{l,m}$ used in the main text.  Put
$f_{l,i}=\gamma u_0^{l,i}$ and define
\[
 F_l:\mathbb C^m\to X,\qquad F_l c=\sum_{i=1}^m c_i f_{l,i}.
\]
The boundary transfer operator has the pole decomposition
\[
 \mathcal N(z)=-\frac{F_lF_l^*}{z-\lambda_0^l}+\mathcal N_l(z),
\]
where the regular part excludes the entire eigenspace.  If
$Q_l=I+\varepsilon\mathcal N_l(z)$, the finite-rank determinant identity
reduces the characteristic equation to
\begin{equation}\label{eq:multiple-effective-matrix}
 \det\left((z-\lambda_0^l)I_m
            -\varepsilon F_l^*Q_l(z,\varepsilon)^{-1}F_l\right)=0.
\end{equation}
Here $F_l^*F_l=C$, with exactly the entries
$C_{ij}=\int_{\partial\Omega}u_0^{l,i}u_0^{l,j}\,d\sigma$ from the main text.
For $|z-\lambda_0^l|=\mathcal O(\varepsilon)$,
\[
 F_l^*Q_l(z,\varepsilon)^{-1}F_l=C+\mathcal O(\varepsilon)
\]
in finite-dimensional operator norm.  All characteristic values in this
cluster satisfy $|z-\lambda_0^l|\leq C_l|\varepsilon|$, as follows directly
from \eqref{eq:multiple-effective-matrix}.  For real $\varepsilon>0$, set
$z=\lambda_0^l+\varepsilon\eta$.  The equation becomes
$\det(\eta I_m-C+\mathcal O(\varepsilon))=0$.
Since $C$ is Hermitian, its resolvent bound shows that every such $\eta$
is within $C_l\varepsilon$ of $\sigma(C)$.  Rouch\'e's theorem on circles
around the distinct eigenvalues of $C$ preserves their multiplicities.
Therefore, if $\mu_1,\ldots,\mu_m$ are the eigenvalues of $C$ counted with
their multiplicities,
\[
 \lambda_\varepsilon^{l,i}
       =\lambda_0^l+\varepsilon\mu_i+\mathcal O(\varepsilon^2),
       \qquad i=1,\ldots,m.
\]
This proves \eqref{asymp:eigenvaluem}.  The simple-mode eigenfunction formula
cannot be applied to an arbitrary vector in a multiple eigenspace.
The combinations determined by $C$ are first chosen; unresolved
multiplicities require the next effective matrix.  No denominator
$\lambda_0^l-\lambda_0^j$ is used within the original cluster.

\subsection{A quantitative remainder controlled by the spectral gap}

The preceding remainder estimates hold for each fixed eigenvalue. To control the remainder as the mode index varies and in the high-frequency discussion, we now make its dependence on the spectral gap explicit.
For a simple eigenvalue, retain $\gamma^l$ from \eqref{def:gammal} and set
\begin{equation}\label{eq:Hl-A}
 H_l:=1+\frac{1+\lambda_0^l}{\gamma^l}.
\end{equation}
Constants in the following proposition depend on the fixed smooth domain,
but not on the index $l$.

\begin{prop}[Gap-dependent Taylor remainder]\label{prop:gap-remainder-A}
There are $c_\Omega,C_\Omega>0$ such that every simple Neumann eigenvalue
satisfies
\begin{equation}\label{eq:uniform-Taylor-A}
 \left|\lambda_\varepsilon^l-\lambda_0^l
                   -\varepsilon a_l-\varepsilon^2B_l\right|
 \leq C_\Omega a_l\varepsilon^3 H_l^2,
 \qquad 0\leq\varepsilon H_l\leq c_\Omega.
\end{equation}
Moreover, $|B_l|\leq C_\Omega a_l H_l$.
\end{prop}
\begin{proof}
For $|z-\lambda_0^l|\leq\gamma^l/4$, the spectral representation
\eqref{eq:N-factorization-A}, with the $l$th mode removed, gives
\begin{equation}\label{eq:Nl-gap-A}
 \|\mathcal N_l(z)\|_{X\to X}
 \leq\|\mathcal T\|^2
      \sup_{j\ne l}\frac{1+\lambda_0^j}{|\lambda_0^j-z|}
 \leq C_\Omega H_l.
\end{equation}
The last inequality follows from
$|\lambda_0^j-z|\geq\tfrac34|\lambda_0^j-\lambda_0^l|$ and
\[
 \frac{1+\lambda_0^j}{|\lambda_0^j-\lambda_0^l|}
       \leq 1+\frac{1+\lambda_0^l}{\gamma^l}.
\]
The ordinary trace estimate gives $a_l\leq C_\Omega(1+\lambda_0^l)$.
Hence, one can choose a domain-dependent constant $c>0$ so that, for all
complex $|\varepsilon|\leq\rho_l:=c/H_l$, both
\[
 \|\varepsilon\mathcal N_l(z)\|\leq\tfrac12,
 \qquad 2|\varepsilon|a_l<\gamma^l/4
\]
hold on the disk.  In obtaining the second inequality, use
$(1+\lambda_0^l)/H_l\leq\gamma^l$.
Then $Q_l^{-1}$ has norm at most $2$ and
\[
 |\varepsilon(f_l,Q_l(z,\varepsilon)^{-1}f_l)|
                          \leq2|\varepsilon|a_l.
\]
Rouch\'e's theorem applied to the Schur complement of
\eqref{eq:regularized-pencil} yields one holomorphic characteristic value
for $|\varepsilon|<\rho_l$, and
$|\lambda_\varepsilon^l-\lambda_0^l|\leq2a_l|\varepsilon|$.
Thus, the removable quotient
\[
 g_l(\varepsilon):=
       \frac{\lambda_\varepsilon^l-\lambda_0^l}{\varepsilon}
\]
is holomorphic and bounded by $2a_l$ on that parameter disk.
Cauchy's coefficient estimates, applied to
$g_l=a_l+B_l\varepsilon+\cdots$, imply
$|B_l|\leq2a_l/\rho_l$ and, for $|\varepsilon|\leq\rho_l/2$,
\[
 |g_l(\varepsilon)-a_l-B_l\varepsilon|
 \leq 4a_l|\varepsilon|^2/\rho_l^2.
\]
Multiplying by $|\varepsilon|$ proves \eqref{eq:uniform-Taylor-A}, after
adjusting the domain-dependent constants.  This proof controls the full
analytic remainder, not only its second coefficient.
\end{proof}

In particular, \eqref{eq:uniform-Taylor-A} implies that
\[
 \lambda_\varepsilon^l
 =\lambda_0^l+\varepsilon a_l+\varepsilon^2B_l
      +\mathcal O_\Omega\!\left(a_l\varepsilon^3
               \left(1+\frac{1+\lambda_0^l}{\gamma^l}\right)^2\right).
\]
Together with $a_l\leq C_\Omega(1+\lambda_0^l)^{1/3}$ for smooth domains, this gives a quantitative interpretation of
\eqref{high-frequency}.  

\bibliographystyle{plain}
\bibliography{ref}
\end{document}